\documentclass[journal]{IEEEtran} % draftcls,

\usepackage[cmex10]{amsmath}
\usepackage{amsfonts,amssymb}
\usepackage{color}
\usepackage{enumerate}
\usepackage{graphicx}
\usepackage{multirow}
\usepackage{placeins}
\usepackage{psfrag}
\usepackage{url}
\usepackage{xspace}
\usepackage{RTMtools}
\allowdisplaybreaks[4] % 4 most permissive breaking of equations across columns/pages

\newcommand{\Equation}{} %{Eqn.} % change style for document class change here
\newcommand{\Equations}{} %{Eqns.} % change style for document class change here
\newcommand{\Figure}{Fig.} % change style for document class change here
\newcommand{\Figures}{Figs.} % change style for document class change here
\newcommand{\pref}[1]{\eqref{#1}}
\newcommand{\Section}{Section}
\newcommand{\ue}{\tilde{u}}
\newcommand{\up}{\hat{u}}
\newcommand{\ve}{\tilde{v}}
\newcommand{\vp}{\hat{v}}
\newcommand{\xe}{\tilde{x}}
\newcommand{\tpo}{t_0^p}
\newcommand{\tpf}{t_f^p}

\newcommand{\tnormSoln}{\overline{t}_{Soln}}
\newcommand{\ml}{\zeta} % mode label
\newcommand{\ML}{Z}
\newcommand{\MIPT}{MIP~\cite{MulParTooBox2004}\xspace}

\newcommand{\diag}{\mathrm{diag}}
\newcommand{\resizeeqnapp}{\small}
\newcommand{\resizematapp}[1]{\small{#1}}
\newcommand{\smat}[1]{\left[\begin{smallmatrix} #1 \end{smallmatrix}\right]}

\begin{document}
\title{A Comparison of the Embedding Method to Multi-Parametric Programming, Mixed-Integer Programming,
Gradient-Descent, and Hybrid Minimum Principle Based Methods}
\author{Richard~T.~Meyer,~\IEEEmembership{Member,~IEEE,}
  Milo\v{s}~\v{Z}efran,~\IEEEmembership{Member,~IEEE,}

  and~Raymond~A.~DeCarlo,~\IEEEmembership{Fellow,~IEEE}% <-this % stops a space
%\thanks{A preliminary version of this paper has been submitted to the
%  2012 IEEE Conference on Decision and Control \cite{prelim-cdc12}.}%
%\presented at the 2007 IEEE Int. Conference on Robotics and Automation{Manuscript received xx}% <-this % stops a space
\thanks{R.T. Meyer is with the School of Electrical and Computer Engineering, Purdue University, West Lafayette,
Indiana 47907 (email: rtmeyer@purdue.edu).}% <-this % stops a space
\thanks{M. \v{Z}efran is with the Department of Electrical and Computer Engineering,
University of Illinois at Chicago, Chicago, Illinois, 60607. His
research was supported in part by NSF grants
    IIS-0905593, CNS-0910988 and CNS-1035914.}% <-this % stops a space
\thanks{R.A. DeCarlo is with the School of Electrical and Computer Engineering
Purdue University, West Lafayette, Indiana 47907.}} %<- double curly brackets to match \author{

% running page headers
%\markboth{IEEE Transactions on Control Systems Technology, Vol.~X No.~X, Month~2012}{Meyer
%\MakeLowercase{\textit{et al.}}: Title chosen} % official header
\markboth{IEEE Transactions on Control Systems Technology}{Meyer
\MakeLowercase{\textit{et al.}}: Title chosen} % header added for RAD

\maketitle
%=========================================================================
\begin{abstract}
In recent years, the embedding approach for solving switched optimal control problems has been developed in a series of papers.  However, the embedding approach, which advantageously converts the hybrid optimal control problem to a classical nonlinear optimization, has not been extensively compared to alternative solution approaches.  The goal of this paper is thus to compare the embedding approach to multi-parametric programming, mixed-integer programming (e.g., CPLEX), and gradient-descent based methods in the context of five recently published examples: a spring-mass system,  moving-target tracking for a mobile robot, two-tank filling, DC-DC boost converter, and skid-steered vehicle. A sixth example, an autonomous switched 11-region linear system, is used to compare a  hybrid minimum principle method and traditional numerical  programming.  For a given performance index for each case, cost and solution times are presented. It is shown that there are numerical  advantages of the embedding approach: lower performance index cost (except in some instances when autonomous switches are present),  generally faster solution time, and convergence to a solution when other methods may fail. In addition,  the embedding method requires no ad~hoc assumptions (e.g.,  predetermined mode sequences) or specialized control  models. Theoretical advantages of the embedding approach over the other methods are also described: guaranteed existence of a solution under mild conditions, convexity of the embedded hybrid optimization problem (under the customary conditions on the performance index), solvability with traditional techniques (e.g., sequential quadratic  programming) avoiding the combinatorial complexity in the number of  modes/discrete variables of mixed-integer programming, applicability  to affine nonlinear systems, and no need to explicitly assign discrete/mode variables to  autonomous switches. Finally,  common misconceptions regarding the embedding approach are  addressed including whether it uses an average value control model  (no), whether it is necessary to ``tweak'' the algorithm to obtain  bang-bang solutions (no), whether it requires infinite switching to  implement embedded solution (no), and whether it has real-time  capability (yes).

\end{abstract}

%=========================================================================
\section{Introduction}

Hybrid and switched systems have modes of operation.  Switches can be controlled or autonomous. In much of the literature, all such modes are labeled with discrete variables.  Solution of these problems has been carried out by a number of optimization algorithms including mixed-integer programming (MIP)~\cite{IntPro1998} and variants of gradient-descent methods~\cite{wardi2010-2137,wardi2012acc,passenberg2010-6666,passenberg2010-4223}. Various implementations of MIP exist, including academic~\cite{MulParTooBox2004} and commercial (CPLEX)~\cite{cplex2012-12d3} solvers. In addition, MIP has been combined with a set of approximations and affine control parameterizations to produce multi-parametric programming (MPP)~\cite{MulParTooBox2004,PreConforLinandHybSys2011}.  An optimization algorithm not evaluated in this work but of importance is relaxed dynamic programming for switched systems~\cite{mariethoz2010-1126,rantzer2006-567,lincoln2006-1249}.
% Response to reviewer 2: The typos have been corrected. Thank you for bringing the availability of a relaxed dynamic programming solver to our attention.  The reference has been added.

%Hybrid and switched systems have modes of operation.  Switches can be controlled or autonomous. In much of the literature, all such modes are labeled with discrete variables.  Solution of these problems is carried out by a number optimization algorithms including academic~\cite{MulParTooBox2004} and commercial (CPLEX)~\cite{cplex2012-12d3} solvers.  mixed-integer programming (MIP) implementations, or variants of gradient-descent methods~\cite{wardi2010-2137,wardi2012acc,passenberg2010-6666,passenberg2010-4223}.  A method not evaluated in this work but of importance is relaxed dynamic programming methods for switched systems~\cite{mariethoz2010-1126,rantzer2006-567,lincoln2006-1249}.

%For simplicity we refer to CPLEX as a separate optimization methodology to distinguish it from MIP implementation in~\cite{MulParTooBox2004}.

In 2005, an embedding approach was developed by Bengea and DeCarlo~\cite{bengea2005-11} and later extended in~\cite{wei2007-264,wei2013-wmr} to a special subclass of autonomous switches.  The embedding approach converts the switched hybrid optimal control problem to a classical continuous optimal control problem. The latter can be converted into a nonlinear programming problem through discretization and solved using sequential-quadratic programming (SQP) for example. This approach has been successfully used in a variety of applications including power management of hybrid electric vehicles~\cite{uthaichana2011-96} and fuel cell hybrid vehicles~\cite{meyer2011-2725}, mobile robot slip control~\cite{wei2007-2373}, and real-time switching control of DC-DC converters~\cite{neely2009-1129,oettmeier2009-3453,neely2010-480,neely2010-compel}.

The goal of this paper is to compare and evaluate the convergence time, resulting performance index costs, and requirements (e.g., requisite assumptions) of these approaches and algorithms in the context of six recently published examples:
\begin{enumerate}[(i)]
\item spring-mass system~\cite{PreConforLinandHybSys2011},
\item a mobile robot~\cite{wardi2010-2137},
\item  two-tank system~\cite{wardi2012acc},
\item a dc-dc boost converter~\cite{mariethoz2010-1126},
\item skid-steered vehicle~\cite{caldwell2010-5423,caldwell2011-50}, and
\item autonomously switched system with 11 state space regions~\cite{passenberg2010-6666,passenberg2010-4223}.
\end{enumerate}
Unless otherwise specified, we pair the embedding methodology and MIP solvers with the model-predictive control (MPC) paradigm so they can be compared to MPP. Furthermore, for simplicity we refer to CPLEX as a separate optimization methodology to distinguish it from MIP implementation in~\cite{MulParTooBox2004}. In example (i) we compare the embedding approach to MPP, \MIPT, and CPLEX and then repeat the experiment for a parameter change that more closely resembles a practical example with a reasonable controller bandwidth.  For examples (ii) and (iii) we compare the embedding method to gradient-descent based methods in~\cite{wardi2010-2137,wardi2012acc}, and in example (iv) we compare the embedding approach to MPP, \MIPT, and CPLEX using the available information in~\cite{mariethoz2010-1126}. For example (v) we replicate the results in~\cite{caldwell2010-5423,caldwell2011-50} since their method is an application of the embedding approach. Finally, for example (vi) we show that the minimum principle solution approach in~\cite{passenberg2010-6666} and traditional numerical programming compare favorably since the problem is one of autonomous switches that do not require a mode designation per se as described in~\cite{wei2007-264}.  Although a special case of the embedding method, we note that there are possible advantages of explicitly parameterizing the autonomous switches.

The example comparisons (which are reproducible) delineated in the subsequent sections, show that at least in the context of the examples studied herein, the embedding method is easy to implement, often has faster convergence, achieves lower costs, and avoids the need for simplifying assumptions or approximate models. Some of the attributes of the embedding approach are that it does not require the use of predefined switching sequences, the generation of offline controller maps, or ad~hoc assumptions on continuous time controls such as a single constant control over extended prediction horizons~\cite{mariethoz2010-1126}. An additional attribute is that the embedding method allows for nonlinear affine systems (those linear in the continuous time control) in contrast to~\cite{PreConforLinandHybSys2011,mariethoz2010-1126} that requires a piecewise linear-affine approximation\footnote{In the latest implementation of the MPP method~\cite{MulParTooBox2004} the authors caution against using general nonlinear models.}.  Finally, for affine systems the embedding approach theoretically guarantees existence of a solution under appropriate convexity of the integrand of the performance index.

%  To the reviewer:  We have rewritten the final paragraph of the introduction and we restrict the comments to the context of the studied examples.  Of course, each reviewer and researcher has their own favorite method and are biased toward that method.  The point of this paper was to address misconceptions in regards to the embedding method that have appeared as biases against the embedding method.  The language of the original paragraph was perhaps too strong, so we have removed the possibly offending adjectives.  Nevertheless, we believe it is important to identify the attributes of the embedding method in comparison to the other methods, which is the thrust of the paper.

%+++++++++++++++++++++++++++++++++++++++++++++++++++++++++++++++++++++++++
\section{Embedding Approach Description}
Details of the embedding method have been set forth in~\cite{bengea2005-11,wei2007-264,uthaichana2011-96,meyer2011-2725,bengea2003-11,bengea2011-311,bengea2003-5925}. This subsection describes the relationship of the embedded and switched system models and that of
the associated switched optimal control problems.

The embedded model of a (possibly nonlinear) switched system takes the form
\begin{equation}\label{e:SwiSysDyn}
  \dot{x}(t)=\sum\limits_{\ml=1}^{d_v}{{{v}_{\ml}}(t){{f}_{\ml}}
    \left(x(t),{{u}_{\ml}}(t) \right)},\quad x({{t}_{0}})={{x}_{0}}
\end{equation}
subject to % could use intertext command here if want to align
\begin{equation}\label{e:EmbModSwiSumOne}
  \sum\limits_{\ml=1}^{d_v}{{{v}_{\ml}}(t)}=1
\end{equation}
where (i) $x(t)\in {{\Real}^{n}}$ is the continuous state; (ii) ${u}_{\ml}(t)\in{{\Real}^{m}}$ is the mode-specific continuous control input that is assumed to be a measurable function; (iii) ${{v}_{\ml}}\in [0,1]\,(\ml=1,...,{d_v})$ is a measurable function that controls mode activation according to \Equation\pref{e:EmbModSwiSumOne}\footnote{Note that if $v_\ml(t)\in\{0,1\}$ and $u(t)=u_\ml(t)$ for all $\ml$, the original switched system is recovered since \Equation\pref{e:EmbModSwiSumOne} forces precisely one mode to be selected.  The key to the usefulness of the embedded formulation is that the trajectories of the switched system ($v_\ml(t)\in\{0,1\}$ and $u(t)=u_\ml(t)$) are dense in the trajectories of the embedded system~\cite{bengea2005-11}.}; and (iv) ${{f}_{\ml}}(x,u):{{\Real}^{n}}\times {{\Real}^{m}}\to {{\Real}^{n}}$ is a piecewise $C^{1}$ vector field where the potential discontinuity accounts for a special subclass of autonomous switches.  See Appendix~\ref{a:DynVecFie} for details. Using the embedded formulation given above with ${v_\ml\in [0,1]\,(\ml=1,\ldots,d_v)}$, the embedded performance index (PI) is defined as
\begin{equation}
\begin{split}
J({{t}_{0}},t_f,{{x}_{0}},& {{u}_{\ml}}\ldots,{{v}_{\ml}}\ldots)=g({{t}_{0}},{{x}_{0}},{{t}_{f}},x({{t}_{f}})) \\ &+\int\limits_{{{t}_{0}}}^{{{t}_{f}}}{\sum\limits_{\ml=0}^{d_v}{{{v}_{\ml}}(t)F_{\ml}^{0}\left( x(t),u_{\ml}(t)\right)dt}}
\end{split}
\end{equation}
This leads to the embedded optimal control problem (EOCP):
\begin{equation}
\underset{u_\ml\in\Omega,v_\ml\in[0,1],\ml\in \{1,\ldots,d_v\}}{\min}\;J({{t}_{0}},t_f,{{x}_{0}},{{u}_{\ml}}\ldots,{{v}_{\ml}}\ldots)
\end{equation}
subject to \Equation\pref{e:SwiSysDyn}-\pref{e:EmbModSwiSumOne}, boundary conditions defined as $\left( {{t}_{0}},x({{t}_{0}}) \right)\in {{T}_{0}}\times {{B}_{0}}$ and $\left({{t}_{f}},x({{t}_{f}}) \right)\in {{T}_{f}}\times {{B}_{f}}$ where ${{T}_{0}}\times {{B}_{0}}\subseteq {{R}^{n+1}}$ and ${{T}_{f}}\times {{B}_{f}}\subseteq {{R}^{n+1}}$ are compact initial and final value sets, ${{u}_{\ml}}(t)\in \Omega$ with $\Omega \subseteq {{\Real}^{m}}$ both compact and convex set, and any other constraints on the controls or state; generally we require that any additional constraints leave the admissible control regions convex.  The EOCP is a classical optimization problem solvable using existing nonlinear programming algorithms; this avoids the search and optimize iterations (combinatorial complexity) associated with mixed-integer programming approaches.  Finding bang-bang solutions, when they exist, is generally not problematic because non-bang-bang solutions require that two or more of the largest mode-specific Hamiltonians be equal as described in [2], which is uncommon.  Typically then, the numerical algorithm will naturally converge to a solution that is bang-bang almost everywhere.  To guarantee existence of a solution to the EOCP, we additionally assume that the set of vector fields are affine in the continuous control variable, i.e.,
\begin{equation}
{{f}_{\ml}}(x,{{u}_{0}})={{A}_{\ml}}(x)+{{B}_{\ml}}(x)\cdot{{u}_{\ml}}
\end{equation}
and the integrands, $F_{\ml}^{0}\left( x,{{u}_{\ml}} \right)$, of the performance index are convex in ${{u}_{\ml}}$ for each $(t,x)$.

%If the switched optimal control problem (SOCP) has a solution, it is a solution of the EOCP except possibly in one circumstance.  That circumstance requires that there be a terminal constraint set and an optimal solution whose value is on the boundary of the terminal constraint set.
If the switched optimal control problem (SOCP), i.e., $v_\ml(t)\in\{0,1\}$ and $u(t)=u_\ml(t)$ for all $\ml$, has a solution, it is a solution of the EOCP except possibly in the following isolated case: there is a terminal constraint set and an optimal solution reaches the terminal constraint set at a point on its boundary. In this rare and easily fixed case, it is possible (but not necessarily true) that the SOCP solution is bounded away from the solution to the EOCP; in this case, the SOCP cost is bounded above the EOCP cost. The ``fix'' is to simply cover the terminal state constraint set with a slightly larger open set; this would result in a ``new'' SOCP solution with cost equal to the EOCP solution. The only example in this paper with a terminal constraint set is in \Section~\ref{s:springmass}, a spring mass system, and we construct a switched solution via the embedding approach meeting the terminal conditions with a cost lower than MPP as published in~\cite{PreConforLinandHybSys2011} and MIP solutions from the multi-parametric toolbox (MPT)~\cite{MulParTooBox2004} and CPLEX~\cite{cplex2012-12d3}.

Thus for the models studied in the examples of this paper, the EOCP always has a solution and if the associated SOCP has a solution, the SOCP solution is one of the possibly non-unique solutions given by the optimal solution space of the EOCP.  The SOCP may not have a solution in the sense that there is no minimum of the performance index over the set $\{0,1\}$ due for example to constraints; in this case, the EOCP solution is the infimum over the set $\{0,1\}$. (Thus the MPP and MIP approaches are infeasible when there does not exist a solution.) Only when there does not exist a SOCP solution would one approximate the EOCP solution with a switched solution using finite-time switching; this is always possible since the switched system trajectories are \textit{dense} in the embedded system trajectories.  Bengea and DeCarlo~\cite{bengea2005-11} provide a construction for approximating the EOCP solution with an SOCP trajectory to any given precision using the Chattering Lemma; practically speaking the duty cycle interpretation is often adequate.  A brief study of projection methods is set forth in~\cite{meyer2013-1}.  All duty cycle interpretations only require finite-time switching. Finally, we point out that when the SOCP does NOT have a solution, i.e., the performance index does not have a minimum over the class of switched systems, the costs associated with the embedded solution and the projection of that solution via a duty cycle interpretation onto the feasible set $\{0,1\}$, (PWM solution) are virtually identical~\cite{neely2009-1129,neely2010-480,neely2010-compel}.

In general, the SOCP is \textit{not} convex.  The EOCP is convex for the reasonable conditions set forth in Bengea and DeCarlo~\cite{bengea2005-11} which are satisfied by the formulations in this paper.  If the SOCP does not have a solution as described previously, then mixed-integer programming approaches are ill-posed. On the other hand, since the EOCP is convex and there are no integer variables (even in the presence of autonomous switches) and since it always has a solution under mild conditions, the EOCP can be solved using classical nonlinear programming techniques such as sequential quadratic programming (SQP). One of the observations of this paper is that MIP methods fail even for some of the simple examples tested here due to combinatorial complexity; instead the embedding approach always finds a solution and is typically faster. The commercially available CPLEX favorably compares to the embedding approach in speed if the number of integer variables is not too large; when there are autonomous switches it sometimes also achieves slightly lower cost.

Finally, we point out to the reader that our implementation of the embedding method uses the Optimization Toolbox (the \textit{fmincon} function) in MATLAB.  The Optimization Toolbox is widely available and has a reasonably good engine, although other engines are known to be better. However, the MATLAB implementation has a considerable I/O overhead\footnote{Our currently unpublished control work on a naval power system indicates that if the overhead is removed,   \textit{fmincon} solution times can be reduced up to an order of magnitude.}. Thus our results do not reflect the full power and achievable speed of the embedding method; but they do provide a reasonable reference for future evaluations.  As such, the experiments should be easily reproducible using commonly available software.

%+++++++++++++++++++++++++++++++++++++++++++++++++++++++++++++++++++++++++
\subsection{Common Concerns about the Embedding Approach}
One common misconception is that the embedding method averages the vector fields of the switched system similar to the Filippov method in variable structure control (VSC).  In VSC, Filippov's method is used to determine a solution to a differential equation whose right hand side is discontinuous on a sliding manifold or discontinuity surface due to infinitely fast switching in the control. Filippov's method forms a convex combination of two vector fields $(1-\alpha)f^+ + \alpha f^-$ and chooses the variable $\alpha$ to achieve an ``average'' value consistent with a tangent plane to the discontinuity surface~\cite{matthews1998-187}. This is not the case in the embedded optimal model where the equivalent of $\alpha$ is a control variable to be chosen so that a performance metric is minimized.

Similarly, the boost and buck converter literature considers time scale separation and linearization about an operating point(s) to obtain an ``average value model''. On the other hand, the embedding method uses the original model and simply forms a convex combination of the vector fields to create a solution space in which the original problem can be solved.  Nowhere is the model averaged in regards to time scales or operating points~\cite{neely2009-1129,oettmeier2009-3453,neely2010-480,neely2010-compel} as illustrated in \Section~\ref{s:DCDCbooCon}.

Another question concerns guaranteeing that the solutions generated by the optimization algorithm are bang-bang. In general, it is not necessary to tweak the algorithm.  Non-bang-bang solutions, in say a two mode system, require that the Hamiltonians in each mode of operation be equal numerically \cite{bengea2005-11} over a time interval of non-zero measure; such solutions are called \emph{singular} because the equality of the Hamiltonians causes a specific function (within a convex combination of the Hamiltonians) given in~\cite{bengea2005-11} to be identically zero on a time interval of non-zero measure. Equal Hamiltonians and existence of solutions were not considered in, e.g., ~\cite{reidinger1999-2466,reidinger1999-3059}. While this occurs rarely in general, it is more frequent when constraints are imposed on the switching set that must remain convex, as in the work of this paper. In a~3 mode or greater system such as in~\cite{meyer2013-1}, two or more Hamiltonians can be equal, but if one of the many is ``larger'', then the solution is bang-bang.  However, as stated earlier, bang-bang solutions are not necessary since non-bang-bang solutions can always be approximated arbitrarily well with bang-bang solutions.

Another issue arises when $v$ is not in the set $\{0,1\}$. A value of $v$ that minimizes the performance index in the interior $(0,1)$ does NOT mean that the SOCP fails to have a solution or that infinite  switching is required to implement the solution.  In~\cite{bengea2005-11}, the value of $v$ is computed for the example in~\cite{xu2000-2683} which is shown to have an infinite number of bang-bang solutions for $v=0.5$. In that example, $v=0.5$ is shown to mean that one must spend equal amounts of time in each mode. Thus a duty cycle interpretation often used for a projection method is consistent with \textit{known} theoretical properties~\cite{bengea2005-11}.

A comment similar to others discussed above argues that the method allows one to use classical SQP, or that the embedding problem is nothing more than a classical nonlinear optimization. As pointed out above, that is precisely why one wants to use the embedding method, because it completely avoids the combinatorial complexity of mixed-integer programming.  As we will see in the examples to follow, mixed-integer programming is generally slow and often does not converge.

%+++++++++++++++++++++++++++++++++++++++++++++++++++++++++++++++++++++++++
\subsection{A Synopsis of the MPP Approach}\label{s:SynMPPapp}
The MPP approach starts with a (discrete time) piecewise affine model of the system (linear in both the state and the control) defined over a set of polytopic regions\footnote{In our synopsis we assume for clarity that (affine) constraints have been appropriately incorporated when constructing these model regions.} $R_i\subset\Real^n$. If the system model is not of this form, an approximate model must be constructed.  According to~\cite{sontag1981-346,mayne2001-ecc,borrelli2005-1709}, the optimal controller for such a system is piecewise affine over polytopic regions for a linear cost function (1-norm and $\infty$-norm) and to a piecewise affine controller over polytopic regions which are subdivided by quadratic forms for a quadratic cost function (2-norm). Starting in a region $R_0$, the authors compute its Chebychev center $x_0$. Given this Chebychev center and a terminal value such as the origin, MPP proceeds to compute the optimal piecewise affine controllers that minimize either a 1-norm, $\infty$-norm, or 2-norm. The optimal trajectory is found by appropriately piecing together the computed sequence of affine controllers of the form $u(k)=K_1x(k)+K_2$. Subsequently, the MPP algorithm must identify a critical region -- a subset of $R_0$ which is consistent with the previously computed sequence of optimal affine controllers. In the case of a 1-norm and $\infty$-norm, the critical region is a polyhedron, while for a 2-norm it is a subset of a polyhedron bounded by a quadratic surface. In turn, the edges of the critical region are used to further subdivide $R_0$. Within each of these subdivisions, a Chebychev center is computed and the process is then repeated~\cite{bemporad2002-3}.

% Response to reviewer 2: We have incorporated the reviewers language
% into our description of MPP.  We thank the reviewer for clarifying
% this subtlety.

The sequence of affine controllers associated with each critical area for the specific PI, optimization horizon, etc., need to be stored. Then for each initial condition in a critical area, the precomputed affine controllers are applied.  In this sense, MPP can be implemented in real-time provided the conditions associated with the precomputed controllers remain the same.  Any parameter change, terminal condition change, time step change, or horizon change requires complete offline recomputation of all controllers.

% Response to reviewer 1: We have removed the following sentence from
% the paper to reduce bias: This is unlike the embedding approach
% which is able to adapt to a changing control problem in real-time as
% demonstrated in~\cite{neely2009-1129}.

Note that when computing the optimal controller, the optimal mode sequence needs to be computed.  The authors do so by examining all possible mode sequences, leading to combinatorial complexity. In addition, critical regions may overlap; this requires further evaluation to determine the controller resulting in least cost~\cite{borrelli2003-TecRep}. It is worth pointing out that the embedding method could be employed for this purpose, eliminating the combinatorial complexity of the original approach.

% %+++++++++++++++++++++++++++++++++++++++++++++++++++++++++++++++++++++++++
% \subsection{Objectives of this Work}
% The aim of this work is to highlight the difference between the
% embedding approach and four other approaches on a set of examples from
% the literature in terms of performance, convergence, computational
% time and ease of use. When possible, the commercially available tool
% CPLEX~\cite{cplex2012-12d3} is also included in the
% comparisons. The experiments suggest that the embedding approach
% almost uniformly outperforms all of the other examined
% approaches in regards to these metrics.

%=========================================================================
\section{Hybrid Optimal Control Examples}
% note to self: MPT can only handle polynomial type nonlinearities in constraints!

In this section, the embedding approach for hybrid optimal control is applied to a spring-mass hybrid system~\cite{PreConforLinandHybSys2011}, a switched-mode mobile robot~\cite{wardi2010-2137}, a two-tank hybrid system~\cite{wardi2012acc}, a DC-DC boost converter~\cite{mariethoz2010-1126}, and a skid-steered vehicle~\cite{caldwell2010-5423,caldwell2011-50}. The results from the embedding approach are compared to those from MPP~\cite{PreConforLinandHybSys2011}, MIP implementation in~\cite{MulParTooBox2004}, CPLEX, the method that computes switching times for a pre-determined mode sequence proposed in~\cite{wardi2010-2137} as well as a follow-on approach that also finds the mode sequence~\cite{wardi2012acc}. Finally, as set forth in~\cite{wei2007-264}, traditional numerical methods are applied to an 11 autonomous mode linear system~\cite{passenberg2010-6666,passenberg2010-4223}. Passenberg\etal apply a hybrid minimum principle~\cite{passenberg2010-6666,passenberg2010-4223} to the problem and by formally labeling autonomous modes to explicitly account for transitions across discontinuity surfaces.

In this work, all the approaches were implemented in MATLAB (version 2010b), except for CPLEX, which is directly called from MATLAB. Appendix~\ref{a:EOCMATSolAlg} outlines the MATLAB-based embedding approach. MPP approach is described in \Section~\ref{s:SynMPPapp}. Each of the other methods is briefly presented in the example where it first appears.

In the following examples, the term ``mode'' indicates a dynamical vector field selected with a discrete control input.  However, ``mode'' has been used in the past to indicate both controlled and uncontrolled (autonomous) switching of vector fields.  Herein, we term the use of ``mode'' resulting from an autonomous switch as an ``autonomous-mode'' or ``a-mode''.  In~\cite{wei2007-264,bengea2011-311}, it was shown that including autonomous switches in the mode definitions is unnecessary for the embedding approach and control problems with only autonomous switches need no mode designations and are solvable with traditional numerical programming.

Further, in each of the examples we use the terms ``numerical optimization cost'' and the ``simulation cost''.  By numerical optimization cost we mean the cost computed via the numerical optimization program using collocation and trapezoidal numerical integration of the PI; by simulated cost we mean the cost obtained by numerically integrating the system ODE's using the piecewise constant continuous controls from the numerical optimization.

In all of the examples, the EOCP is solved using MATLAB's \textit{fmincon} function following the general procedure outlined in Appendix~\ref{a:EOCMATSolAlg}. Using a numerically superior solver would only improve the convergence rate and solution times of the results reported for the EOCP.  The MPP and \MIPT results were obtained using version~2.63 of the MPT\footnote{Improvements to MPP have been made to reduce computational complexity~\cite{tondel2003-945,jones2010-2542}.  However, we restricted ourselves to using the MPT commonly available to the engineering community.}.  Also, CPLEX version~12.4 was used to obtain CPLEX results.  All the code used to generate the results in this paper can be accessed at~\cite{CodeURL}.

%+++++++++++++++++++++++++++++++++++++++++++++++++++++++++++++++++++++++++
\subsection{Spring-Mass Hybrid System~\cite{PreConforLinandHybSys2011}}\label{s:springmass}
Example~14.2 in~\cite{PreConforLinandHybSys2011} introduces the spring-mass hybrid optimal control problem.  A mass is connected to ground with a spring in series with a damper that represents viscous friction. The spring has affine characteristics and the viscous friction coefficient can be changed from one value $b_1$ to a different value $b_2$ instantaneously with a binary input. The continuous-time spring-mass system dynamics are
\begin{align}
\begin{split}\label{e:springmass-dyn-ct}
\dot{x}_1(t)=&x_2(t)\\
M\dot{x}_2(t)=&-k(x_1(t))-b(u_2(t))x_2(t)+u_1(t)
\end{split}\\
\intertext{where the the spring coefficient $k(x_1(t))$ and viscous friction coefficient $b(u_2(t))$ are}
k(x_1(t))=&\begin{cases}
    k_1x_1(t)+d_1,&\; x_1(t)\leq x_m\\
    k_2x_1(t)+d_2,&\; x_1(t)>x_m
    \end{cases}\label{e:springmass-k}\\
b(u_2(t))=&\begin{cases}
    b_1,&\; u_2(t)=1\\
    b_2,&\; u_2(t)=0
    \end{cases}
\end{align}
with $x_1$ and $x_2$ the mass position and velocity, respectively; the system has two modes and two a-modes.  System parameters are given as $M=1$, $b_1=1$, $b_2=50$, $k_1=1$, $k_2=3$, $d_1=1$, $d_2=7.5$, and $x_m=1$.

The spring-mass example is solved in the context of both MPC and a single optimization. Within the MPC context, the prediction horizon has $N$ sample intervals\footnote{The end-point conditions are treated using a shrinking horizon, except for MPP where this is not easily done.}. In each case, the optimizations are appropriately reformulated and solved with (i) MPP\footnote{The MPP implementation  relies on the MIP solver provided in~\cite{MulParTooBox2004}.} (the original solution approach in~\cite{PreConforLinandHybSys2011} is for MPC); (ii) MIP as implemented in~\cite{MulParTooBox2004}; (iii) CPLEX; and (iv) the embedding approach.

To implement the MPP algorithm (described in \Section~\ref{s:SynMPPapp}), four ``modes'' were defined depending on the discrete input $u_2$ and mass position, $x_1$, which determines the spring coefficient with \Equation\pref{e:springmass-k}; the ``modes'' are a mixture of both modes in the sense used here and a-modes. The dynamics in each ``mode'' were discretized with a $0.5$ time unit sampling resulting in a discrete-time PWA system:
\begin{equation}\label{e:springmass-dyn-dt}
\begin{split}
&x(k+1)\\
&=\begin{cases}
    A_1x(k)+B_1u_1(k)+F_1,& x_1(k)\leq 1,u_2(k)\leq0.5\\
    A_2x(k)+B_2u_1(k)+F_2,& x_1(k)> 1,u_2(k)\leq0.5\\
    A_3x(k)+B_3u_1(k)+F_3,& x_1(k)\leq 1,u_2(k)\geq0.5\\
    A_4x(k)+B_4u_1(k)+F_4,& x_1(k)> 1,u_2(k)\geq0.5
\end{cases}
\end{split}
\end{equation}
where $A_i$, $B_i$, and $F_i$ are listed in
Appendix~\ref{a:SprMasDisTimMod}.  The MPP control objective is to minimize 
% Response to reveiwer 2: We have removed the footnote regarding the matrices and now reference the book as being forthcoming.  Please note that we notified the book authors.  
\begin{equation}\label{e:springmass-J-dt-mpp}
\begin{split}
J_m(x(0),U_0,&N)=x(N)^TPx(N) \\
&+\sum_{k=0}^{N-1}x(k)^TQx(k)+u(k)^TRu(k)
\end{split}
\end{equation}
over $U_0=[u(0)^T,\ldots,u(N-1)^T]^T$ ($u(k)=[u_1(k),u_2(k)]^T$) and
subject to \Equation\pref{e:springmass-dyn-dt}, $x_1,x_2\in[-5,5]$,
$u_1\in[-10,10]$, and a terminal state constraint set
$\mathcal{X}_f\in[-0.01,0.01]\times[-0.01,0.01]$.  The original
performance index parameters are $N=3$, $P=Q=I_2$ and
$R=\left[\begin{smallmatrix} 0.2 & 0\\0 & 1\end{smallmatrix}\right]$
and the initial state is $x(0)=[3,4]^T$.  For the simulations here, $R$ is changed to $\left[\begin{smallmatrix}0.2 & 0\\0 &  0\end{smallmatrix}\right]$ because the current embedding approach theory does not consider penalized mode switch values in the PI.

For the alternate embedded model, there are only two (controlled) modes, since the embedding methodology does not require parameterization of autonomous switching~\cite{wei2007-264} (a-modes).
\begin{equation}\label{e:springmass-dyn-ct-eocp}
\begin{split}
\dot{\xe}(t)=&(1-\ve(t))\bmat{\xe_2(t)\\ -\frac{1}{M}k(\xe_1)-\frac{b_2}{M}\xe_2(t)+\frac{1}{M}\ue_1^0(t)}\\
&+\ve(t)\bmat{\xe_2(t)\\ -\frac{1}{M}k(\xe_1)-\frac{b_1}{M}\xe_2(t)+\frac{1}{M}\ue_1^1(t)}
\end{split}
\end{equation}
where $(\tilde{\cdot})$ is an embedded system value, $\ve\in[0,1]$ is the embedded mode switch value and $\ue_1^0$ and $\ue_1^1$ are the continuous controls in each mode.  Notice that $\ve=0$/$\ve=1$ in \Equation\pref{e:springmass-dyn-ct-eocp} corresponds to $u_2=0$/$u_2=1$ in \Equation\pref{e:springmass-dyn-ct}.  The embedded problem is
\begin{gather}
\min_{\ue_1^0,\ue_1^1,\ve}\; J_E(x(\tpo),\ue_1^0,\ue_1^1,\ve,\tpo,\tpf)\\
\intertext{where $[\tpo,\tpf]$ is the prediction interval,}
\begin{split}\label{e:springmass-J-ct-eocp}
J_E=\int_{\tpo}^{\tpf}&\left\{\xe(t)^TQ_E(t)\xe(t)\right.\\
&+(1-\ve(t))R_{E11}(t)(\ue_1^0(t))^2\\
&\left.+\ve(t)R_{E11}(t)(\ue_1^1(t))^2\right\}dt,
\end{split}
\end{gather}
and the minimization is subject to \Equation\pref{e:springmass-dyn-ct-eocp}, the previously defined bounds on the states and continuous control input, and the terminal state constraint set. The numerical optimization requires a discretized embedded system representation of \Equation\pref{e:springmass-dyn-dt}:
\begin{equation}
\begin{split}
&\xe(k+1)\\
&=\begin{cases}
    (1-\ve(k))\left[A_1\xe(k)+B_1\ue_1^0(k)+f_1\right]\\
    +\ve(k)\left[A_3\xe(k)+B_3\ue_1^1(k)+f_3\right],& \; \xe_1(k)\leq 1\\
    (1-\ve(k))\left[A_2\xe(k)+B_2\ue_1^0(k)+f_2\right] \\
    +\ve(k)\left[A_4\xe(k)+B_4\ue_1^1(k)+f_4\right], & \; \xe_1(k)> 1
\end{cases}
\end{split}
\end{equation}
The trapezoidal numerical integration of \Equation\pref{e:springmass-J-ct-eocp} is
\begin{equation}\label{e:springmass-J-dt-eocp}
\begin{split}
\hat{J}_E=&\xe(N)^TP\xe(N) +\sum_{k=0}^{N-1}\left\{\xe(k)^TQ\xe(k)\right.\\
&+(1-\ve(k))R_{11}(\ue_1^0(k))^2+\ve(k)R_{11}(\ue_1^1(k))^2
\end{split}
\end{equation}
where $Q_E(t-\tpo)$ and $R_{E11}(t-\tpo)$ are chosen such that \Equation\pref{e:springmass-J-dt-eocp} is equivalent to \Equation\pref{e:springmass-J-dt-mpp} and $R_{11}$ represents the value of $R(1,1)$. If $\ve(k)\in(0,1)$, then mode and control projections are required. Here, the projection method of~\cite{meyer2011-2725,meyer2013-1} is applied where the projected mode, $\vp$, is zero if $\ve\leq0.5$ and one otherwise.  The projected control $\up_1$ is equal to $(1-\ve)\ue_1^0$ if $\vp=0$ and $\ve\ue_1^1$ if $\vp=1$.

\begin{figure}[t]\centering% This file is generated by the MATLAB m-file laprint.m. It can be included
% into LaTeX documents using the packages graphicx, color and psfrag.
% It is accompanied by a postscript file. A sample LaTeX file is:
%    \documentclass{article}\usepackage{graphicx,color,psfrag}
%    \begin{document}\input{f-eoc-mpp-springmass-states}\end{document}
% See http://www.mathworks.de/matlabcentral/fileexchange/loadFile.do?objectId=4638
% for recent versions of laprint.m.
%
% created by:           LaPrint version 3.16 (13.9.2004)
% created on:           18-Feb-2013 13:32:16
% eps bounding box:     10 cm x 7.5 cm
% comment:              
%
\begin{psfrags}%
\psfragscanon%
%
% text strings:
\psfrag{s06}[b][b]{\color[rgb]{0,0,0}\setlength{\tabcolsep}{0pt}\begin{tabular}{c}$x_1(t)$\end{tabular}}%
\psfrag{s07}[b][b]{\color[rgb]{0,0,0}\setlength{\tabcolsep}{0pt}\begin{tabular}{c}$x_2(t)$\end{tabular}}%
\psfrag{s08}[t][t]{\color[rgb]{0,0,0}\setlength{\tabcolsep}{0pt}\begin{tabular}{c}Time (s)\end{tabular}}%
%
% xticklabels:
\psfrag{x01}[t][t]{0}%
\psfrag{x02}[t][t]{2}%
\psfrag{x03}[t][t]{4}%
\psfrag{x04}[t][t]{6}%
\psfrag{x05}[t][t]{8}%
\psfrag{x06}[t][t]{10}%
\psfrag{x07}[t][t]{12}%
\psfrag{x08}[t][t]{0}%
\psfrag{x09}[t][t]{2}%
\psfrag{x10}[t][t]{4}%
\psfrag{x11}[t][t]{6}%
\psfrag{x12}[t][t]{8}%
\psfrag{x13}[t][t]{10}%
\psfrag{x14}[t][t]{12}%
%
% yticklabels:
\psfrag{v01}[r][r]{-5}%
\psfrag{v02}[r][r]{0}%
\psfrag{v03}[r][r]{5}%
\psfrag{v04}[r][r]{-1}%
\psfrag{v05}[r][r]{0}%
\psfrag{v06}[r][r]{1}%
\psfrag{v07}[r][r]{2}%
\psfrag{v08}[r][r]{3}%
%
% Figure:
\resizebox{8cm}{!}{\includegraphics{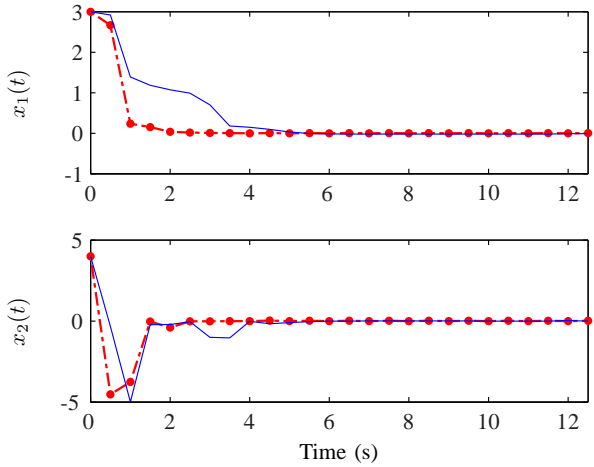}}%
\end{psfrags}%
%
% End f-eoc-mpp-springmass-states.tex
\caption{Spring-mass example state responses:
({\color{Matlabblue}---}) EOC with projection, ({\color{Matlabred}$\bullet$}) MPP (\MIPT and
CPLEX have similar results).}\label{f:eocp-mpp-springmass-states}\end{figure}

\begin{figure}[t]\centering% This file is generated by the MATLAB m-file laprint.m. It can be included
% into LaTeX documents using the packages graphicx, color and psfrag.
% It is accompanied by a postscript file. A sample LaTeX file is:
%    \documentclass{article}\usepackage{graphicx,color,psfrag}
%    \begin{document}\input{f-eoc-mpp-springmass-control}\end{document}
% See http://www.mathworks.de/matlabcentral/fileexchange/loadFile.do?objectId=4638
% for recent versions of laprint.m.
%
% created by:           LaPrint version 3.16 (13.9.2004)
% created on:           18-Feb-2013 13:19:09
% eps bounding box:     10 cm x 7.5 cm
% comment:
%
\begin{psfrags}%
\psfragscanon%
%
% text strings:
\psfrag{s09}[b][b]{\color[rgb]{0,0,0}\setlength{\tabcolsep}{0pt}\begin{tabular}{c}$\hat{v}(t)$, $\tilde{v}(t)$\end{tabular}}%
\psfrag{s10}[b][b]{\color[rgb]{0,0,0}\setlength{\tabcolsep}{0pt}\begin{tabular}{c}$\hat{v}(t)$, $u_2(t)$\end{tabular}}%
\psfrag{s11}[b][b]{\color[rgb]{0,0,0}\setlength{\tabcolsep}{0pt}\begin{tabular}{c}$\hat{u}_1(t)$, $u_1(t)$\end{tabular}}%
\psfrag{s12}[t][t]{\color[rgb]{0,0,0}\setlength{\tabcolsep}{0pt}\begin{tabular}{c}Time (s)\end{tabular}}%
%
% xticklabels:
\psfrag{x01}[t][t]{0}%
\psfrag{x02}[t][t]{2}%
\psfrag{x03}[t][t]{4}%
\psfrag{x04}[t][t]{6}%
\psfrag{x05}[t][t]{8}%
\psfrag{x06}[t][t]{10}%
\psfrag{x07}[t][t]{12}%
\psfrag{x08}[t][t]{0}%
\psfrag{x09}[t][t]{2}%
\psfrag{x10}[t][t]{4}%
\psfrag{x11}[t][t]{6}%
\psfrag{x12}[t][t]{8}%
\psfrag{x13}[t][t]{10}%
\psfrag{x14}[t][t]{12}%
\psfrag{x15}[t][t]{0}%
\psfrag{x16}[t][t]{2}%
\psfrag{x17}[t][t]{4}%
\psfrag{x18}[t][t]{6}%
\psfrag{x19}[t][t]{8}%
\psfrag{x20}[t][t]{10}%
\psfrag{x21}[t][t]{12}%
%
% yticklabels:
\psfrag{v01}[r][r]{0}%
\psfrag{v02}[r][r]{2}%
\psfrag{v03}[r][r]{4}%
\psfrag{v04}[r][r]{6}%
\psfrag{v05}[r][r]{0}%
\psfrag{v06}[r][r]{0.5}%
\psfrag{v07}[r][r]{1}%
\psfrag{v08}[r][r]{0}%
\psfrag{v09}[r][r]{0.5}%
\psfrag{v10}[r][r]{1}%
%
% Figure:
\resizebox{8cm}{!}{\includegraphics{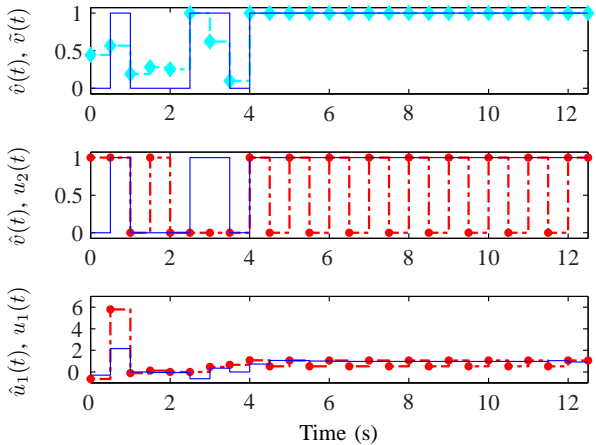}}%
\end{psfrags}%
%
% End f-eoc-mpp-springmass-control.tex
\caption{Spring-mass example continuous and discrete
controls: ({\color{Matlabblue}---}) EOC with projection, ({\color{Matlabcyan}$\blacklozenge$}) EOC
without projection, ({\color{Matlabred}$\bullet$}) MPP (\MIPT and CPLEX have similar
results).}\label{f:eocp-mpp-springmass-control}\end{figure}

\setlength{\tabcolsep}{3.5pt}
\begin{table}[t]
\caption{Spring-mass example EOC (with and without projection), MPP, \MIPT, and CPLEX performance index
cost and solution times for a simulation with $N=3$ and $t_s=0.5$.}
\label{t:SprMasSolCom}
\centering
\begin{tabular}{l c c c c}
\hline
\multirow{2}{*}{Method} & Simulation & Offline & Online & Total\\
 & Cost & $\tnormSoln$ & $\tnormSoln$ & $\tnormSoln$\\
\hline
EOC-Proj. & 71.21 & - & 1.09 & 1.09\\ % normalized 
EOC-No Proj. & 67.44 & - & 1.69 & 1.69\\ % normalized
MPP & 76.43 & 5.69 & 0.27 & 5.96\\ % normalized
\MIPT & 76.20 & 1.10 & 22.62 & 23.72\\ % normalized
CPLEX & 76.43 & - & 1.00 & 1.00\\ % normalized 2.00
\hline
\end{tabular}
\end{table}

\begin{table}[t]
\caption{Spring-mass example EOC (with and without projection, second
  optimization) and CPLEX performance index cost and solution times
  for a simulation with $N=25$ and $t_s=0.5$ (``X'' means the
  algorithm did not converge).}
\label{t:SprMasSolComSinSho}
\centering
\begin{tabular}{l c c}
\hline
\multirow{2}{*}{Method} & Simulation & \multirow{2}{*}{$\tnormSoln$}\\
& Cost & \\
\hline
EOC-Proj. & 105.28 & 1.00 \\ % normalized 30.87
EOC-No Proj. & 53.18 & 1.00 \\ % normalized
EOC-Proj./Opt. & 84.38 & 1.43 \\ % normalized
MPP & X & X \\
\MIPT & X & X \\
CPLEX & 66.92 & 7.01 \\ % normalized
\hline
\end{tabular}
\end{table}

\begin{table}[t]
\caption{Spring-mass example EOC, \MIPT, and CPLEX performance index
 cost and solution times for the redefined system, $N=7$, and
 $t_s=0.2$ (``X'' means the algorithm did not converge).}
\label{t:SprMasSolComRedSys}
\centering
\begin{tabular}{l c c c c}
\hline
\multirow{2}{*}{Method} & Simulation & Offline & Online & Total\\
 & Cost & $\tnormSoln$ & $\tnormSoln$ & $\tnormSoln$ \\
\hline
EOC & 192.04 & - & 1.00 & 1.00 \\ % normalized 12.30
MPP & X & X & X & X \\
\MIPT & 197.23 & 0.21 & 1237.56 & 1237.77 \\ % normalized
CPLEX & 188.09 & - & 1.21 & 1.21 \\
\hline
\end{tabular}
\end{table}

\Figures~\ref{f:eocp-mpp-springmass-states} and~\ref{f:eocp-mpp-springmass-control} shows the state trajectories and evolution of the continuous and discrete controls for the embedding approach and MPP over a 12.5~s simulation with $N=3$ and $t_s=0.5$; \MIPT and CPLEX results are nearly identical to the MPP results and are not shown.  The embedding approach results in~7 mode projections over 25 partitions.  Also, it is seen that the EOC projected mode switch value and MPP mode switch value are different for 11 of the 25 partitions.  Table~\ref{t:SprMasSolCom} lists the PI costs, offline solution times (if applicable), solution times during the system simulation, and total solution times where $\tnormSoln$ indicates time values normalized to the least total solution time observed in the comparison.  For this example, we observe that the embedding approach, both with and without projection, achieves a lower cost than MPP, \MIPT, and CPLEX.  Further, CPLEX solved the problem 9\% faster than the embedding approach with mode projection.

A second test was performed to check solution method performance over longer prediction horizons.  In this test, the control problem was solved once over $N=25$ with $t_s=0.5$, the entire simulation length, and then the system simulated using the resulting controls.  Table~\ref{t:SprMasSolComSinSho} shows the results of this test, again $\tnormSoln$ indicates time values normalized to the least solution time value observed in the comparison.  No MPP or \MIPT results are reported because the MPT software failed to provide results; MPP caused the computer to run out of memory (machines with up to 8~GB of RAM were tried) and \MIPT was stopped after an hour because its screen output indicated no iterations in the first control solution had occurred.  Regardless, CPLEX results in a lower PI cost than either of the embedding approaches.  For this test, after obtaining the embedded approach solution, we took the additional step of solving a second embedded optimization problem with the difference being the mode switch values are set equal to the projected values obtained previously. The new problem solution provides the optimal continuous controls associated with the projected mode sequence.  The embedded approach solution results in 8 mode projections over the 25 partition simulation.  Further, the embedded solution cost without projection is lower than the CPLEX cost which is expected since the embedded solution cost is the infimum.  Also, this test shows the ``curse of dimensionality'' associated with mixed-integer programming methods; the embedded approach solution time is approximately 7 times less than that for CPLEX.

The spring-mass problem as originally presented has inadequate sampling given the bandwidth of the system for the given parameters. Because of an eigenvalue of the continuous time system at approximately $-50$ due to excessively high damping, the sampling interval of $0.5$s does not result in an adequate discrete-time representation of the continuous-time system; a sampling interval of around $0.002$s would be more appropriate.  However, before attempting a sample period of $0.002$, a sample period of $0.05$ with $N=20$ was attempted with the outcome that both MPP and \MIPT did not provide solutions. The embedding approach solution with $15$ projections had a simulation cost of 300.44 and embedded cost of $265.22$.  CPLEX returned a simulation cost of $269.70$ and took $280$ times longer to solve the problem.  Again, the embedded cost was less than CPLEX's cost.
%, and took $380.27$s.  $1.068\dEE{5}$s. CPLEX took over 

In order to obtain a more physically meaningful system, damping values were changed to $b_1=0.1$ and $b_2=0.2$ and the spring constant to $k_2=0.005$.  This results in a minimum time constant of $5.85$.  A sample period of $0.2$ and $N=7$ and $N=7$ provides enough time in the prediction horizon to drive the system states to the terminal constraint set.  Table~\ref{t:SprMasSolComRedSys} lists the performance of the embedding approach, MIP, \MIPT, and CPLEX; the embedding approach solution was bang-bang, making mode projections unnecessary.  The embedding approach solution time is approximately 80\% that of CPLEX and $1240$ times less than that from \MIPT. However, CPLEX results in about a 2\% lower cost than the embedding approach.  It was observed that over several CPLEX solutions it found that one of the optimal values of $x_1$ over $N$ to be at $1$, the autonomous switching point.  Over these same solutions in time, the embedding approach avoided setting $x_1=1$.  If in the embedding solution, the appropriate $x_1$ over $N$ is forced to be $1$, then the embedding approach cost is equal to the CPLEX cost.  In \Section~\ref{s:11StaSpaRegAutSwiExa} and in the Discussion and Conclusions section the effect of autonomous switches on the embedding approach is described in more detail with a suggestion for future research.

% Table~\ref{t:springmassSolutionCompare} lists the performance index
% costs, $J$, and control solution times, $\Delta t^r$, for the
% simulation associated with $N=3$ plus simulations with $N=4$ with
% $t_s=0.375$ and $N=5$ with $t_s=0.3$. For $N=4$ and $N=5$ the total
% simulation times are~12.375 and 12.6, respectively. There are 10
% mode projections over 33 partitions needed for $N=4$ and 7 over 42
% partitions for $N=5$.  We conclude that for the same performance
% metrics and the same example, the embedding approach achieved lower
% costs in all trials. Furthermore, for an MPC window of 5 partitions,
% the embedding approach was 143 times faster than the \MIPT
% approach and 22 times faster than the MPP approach. Although the MPP
% approach generates a look-up table, if any parameters change or
% there are modeling issues, the entire table must be recomputed.  The
% embedding approach was on average about 3.5 times slower than using
% CPLEX. CPLEX is not a general purpose minimization program like
% \textit{fmincon} that we use to solve the EOCP. CPLEX is
% specifically tailored to minimize a quadratic cost function subject
% to linear equality constraints (as is the spring-mass problem).

%+++++++++++++++++++++++++++++++++++++++++++++++++++++++++++++++++++++++++
\subsection{Mobile Robot Hybrid System~\cite{wardi2010-2137}}
Wardi\etal~\cite{wardi2010-2137} consider the problem of a mobile robot tracking a moving target while avoiding obstacles.  The robot has three operating modes: go to goal (G2G), avoid obstacle~1 (Avoid1), and avoid obstacle~2 (Avoid2).  The robot's dynamic equation of motion is
\begin{equation}\label{e:MobRobDyn}
\dot{x}(t)=\bmat{V\cos(x_3(t)) \\ V\sin(x_3(t)) \\u-x_3(t)}
\end{equation}
where $V$ is
\begin{equation}
V=\begin{cases}
    \overline{V},&\; \|x_R-x_G\|\geq r\\
    \frac{\overline{V}}{r}\|x_R-x_G\|,&\;\|x_R-x_G\|<r
\end{cases}
\end{equation}
and $u$ is defined for the three modes as
\begin{align}
u_{G2G}=&\tan^{-1}\left(\frac{x_{G,2}(t)-x_{R,2}(t)}{x_{G,1}(t)-x_{R,1}(t)}\right)\\
u_{Avoid1}=&\begin{cases}
    \phi_{\Phi}-\pi/2,&\; \phi_{\Phi}-x_3\geq 0\\
    \phi_{\Phi}+\pi/2,&\; \phi_{\Phi}-x_3<0
\end{cases}\\
u_{Avoid2}=&\begin{cases}
    \phi_{\Psi}-\pi/2,&\; \phi_{\Psi}-x_3\geq 0 \\
    \phi_{\Psi}+\pi/2,&\; \phi_{\Psi}-x_3<0
\end{cases}
\end{align}
and $\phi_{\Phi}$ and $\phi_{\Phi}$ are
\begin{align}
\phi_{\Phi}=&\tan^{-1}\left(\frac{x_{\Phi,2}-x_{R,2}(t)}{x_{\Phi,1}-x_{R,1}(t)}\right)\\
\phi_{\Psi}=&\tan^{-1}\left(\frac{x_{\Psi,2}-x_{R,2}(t)}{x_{\Psi,1}-x_{R,1}(t)}\right)
\label{e:MobRobPhiPsi}
\end{align}
where $x=[x_R,x_3]^T$ such that $x_R=[x_{R,1},x_{R,2}]^T$ is the global position (two coordinates) of the mobile robot and $x_3\in[0,2\pi)$ is the robot heading angle; $x_G$ is the coordinate pair of the goal; $\overline{V}=1$; $r=0.5$; $x_{\Phi}=[0,4]^T$ are the coordinates of obstacle~1; and $x_{\Psi}=[6,6]^T$ are the coordinates of obstacle~2.

The selection of the robot mode is formulated as an MPC problem with a shrinking horizon in~\cite{wardi2010-2137}. The prediction horizon is initially $[0,t_f]$ divided into $N_0$ sample intervals of length $t_s=t_f/N_0$.  The second horizon is $[t_s,t_f]$ with $N_1=N_0-1$ partitions of length $t_s$.  The pattern of shrinking the horizon continues until it reaches the single partition $[t_f-t_s,t_f]$.  Specifically, the MPC problem is to select G2G, Avoid1, or Avoid2 that minimizes
\begin{equation}\label{e:mobilerobot-J}
\begin{split}
J(x(t),x_G(t),t)=&\int_t^{t_f}\left\{\vphantom{\left(-\frac{\|x_R(t)-x_{\Psi}\|^2}{\beta_12}\right)}
\rho\|x_R(t)-x_G(t)\|^2 \right.\\
&+\alpha_1\exp\left(-\frac{\|x_R(t)-x_{\Phi}\|^2}{\beta_1}\right)\\
&\left.+\alpha_2\exp\left(-\frac{\|x_R(t)-x_{\Psi}\|^2}{\beta_2}\right)\right\}dt
\end{split}
\end{equation}
subject to \Equations\pref{e:MobRobDyn}-\pref{e:MobRobPhiPsi} from the start time $t$ (a multiple of $t_s$) to $t_f$ where $t_f=18$, $N_0=180$, $t_s=0.1$, $\rho=0.1$, $\alpha_1=\alpha_2=500$, and $\beta_1=\beta_2=0.8$. The goal coordinates are not known \textit{a priori}, rather an estimate,
$\tilde{x}_G$ at $s\geq t$ is used in the cost function:
\begin{equation}
\tilde{x}_G(s,t,x_G(t)):=x_G(t)+\dot{x}_G(t)(s-t)
\end{equation}
where the derivative of $x_G(t)$ is approximated with $(x_G(t)-x_G(t-t_s))/t_s$.

Wardi\etal~\cite{wardi2010-2137} solve the above MPC problem by fixing a mode sequence \textit{a priori} and determining the switching times.  After the initial MPC optimization using a gradient-descent technique, if the first mode switching time $t_1\in[0,t_s]$, then it becomes fixed; if $t_1\not\in[0,t_s]$, then its value can still be modified.  The MPC optimization using the gradient-descent technique then continues for shrinking horizons until all switching times are found.  

The embedding approach is applied to the above MPC problem and the results compared to those in~\cite{wardi2010-2137}.  MPP, \MIPT, and CPLEX are not applicable to this problem for two reasons: (i) the dynamics are nonlinear and (ii) the PI contains exponentials of the states.
 
% approach to solving the MPC problem here is to find the switching times for a pre-determined mode sequence.   During a control update, a gradient-descent technique is used to regularly update the mode switching times for upcoming modes left in the sequence.  Their results are reported further on.
%
%The MPC problem in this example is not solved with MPP, \MIPT, and CPLEX because they are not directly applicable to problems like this with nonlinear dynamics and constraints and the exponential function in the PI.  The embedding approach has no such limitations.  Hence only the embedding approach solution to the MPC problem is compared with the published results in~\cite{wardi2010-2137}.

For the embedding model the three robot modes are controlled by $\ve_0$, $\ve_1$, and $\ve_2$:
\begin{equation}\label{e:mobilerobot-dyn-ct-eocp}
\begin{split}
\dot{\xe}(t)=&\ve_0(t)\bmat{V\cos(\xe_3(t)) \\ V\sin(\xe_3(t)) \\u_{G2G}-\xe_3(t)}
+\ve_1(t)\bmat{V\cos(\xe_3(t)) \\ V\sin(\xe_3(t)) \\u_{Avoid1}-\xe_3(t)}\\
&+\ve_2(t)\bmat{V\cos(\xe_3(t)) \\ V\sin(\xe_3(t)) \\u_{Avoid2}-\xe_3(t)}
\end{split}
\end{equation}
where $\ve_i\in[0,1]$ and $\ve_0+\ve_1+\ve_2=1$.  The embedded control problem is to also minimize \Equation\pref{e:mobilerobot-J} over the $\ve_i$ subject to \Equation\pref{e:mobilerobot-dyn-ct-eocp} and the summation constraint on $\ve_i$.  MATLAB's \textit{fmincon} function is used to solve the problem after converting the continuous-time equations to discrete-time equality constraints via collocation~\cite{wei2007-264}.  The cost function is approximated with trapezoidal numerical integration (see Appendix~\ref{a:EOCMATSolAlg} for more detail).  A continuation approach~\cite{richter1983-660} is used to construct an initial guess; furthermore, the continuation approach itself is initialized using Appendix~\ref{a:EOCMATSolAlg}-Step~\ref{el:fminconinitialguess}. When needed, mode projection is accomplished by selecting the embedded mode switch with the largest value as was done in~\cite{meyer2013-1}.

\begin{figure}[]\centering% This file is generated by the MATLAB m-file laprint.m. It can be included
% into LaTeX documents using the packages graphicx, color and psfrag.
% It is accompanied by a postscript file. A sample LaTeX file is:
%    \documentclass{article}\usepackage{graphicx,color,psfrag}
%    \begin{document}\input{f-mobilerobot-trajectory}\end{document}
% See http://www.mathworks.de/matlabcentral/fileexchange/loadFile.do?objectId=4638
% for recent versions of laprint.m.
%
% created by:           LaPrint version 3.16 (13.9.2004)
% created on:           18-Feb-2013 14:00:03
% eps bounding box:     10 cm x 7.5 cm
% comment:              
%
\begin{psfrags}%
\psfragscanon%
%
% text strings:
\psfrag{s03}[b][b]{\color[rgb]{0,0,0}\setlength{\tabcolsep}{0pt}\begin{tabular}{c}$x_2$\end{tabular}}%
\psfrag{s04}[t][t]{\color[rgb]{0,0,0}\setlength{\tabcolsep}{0pt}\begin{tabular}{c}$x_1$\end{tabular}}%
%
% xticklabels:
\psfrag{x01}[t][t]{0}%
\psfrag{x02}[t][t]{2}%
\psfrag{x03}[t][t]{4}%
\psfrag{x04}[t][t]{6}%
\psfrag{x05}[t][t]{8}%
\psfrag{x06}[t][t]{10}%
\psfrag{x07}[t][t]{12}%
%
% yticklabels:
\psfrag{v01}[r][r]{0}%
\psfrag{v02}[r][r]{2}%
\psfrag{v03}[r][r]{4}%
\psfrag{v04}[r][r]{6}%
\psfrag{v05}[r][r]{8}%
\psfrag{v06}[r][r]{10}%
\psfrag{v07}[r][r]{12}%
%
% Figure:
\resizebox{8cm}{!}{\includegraphics{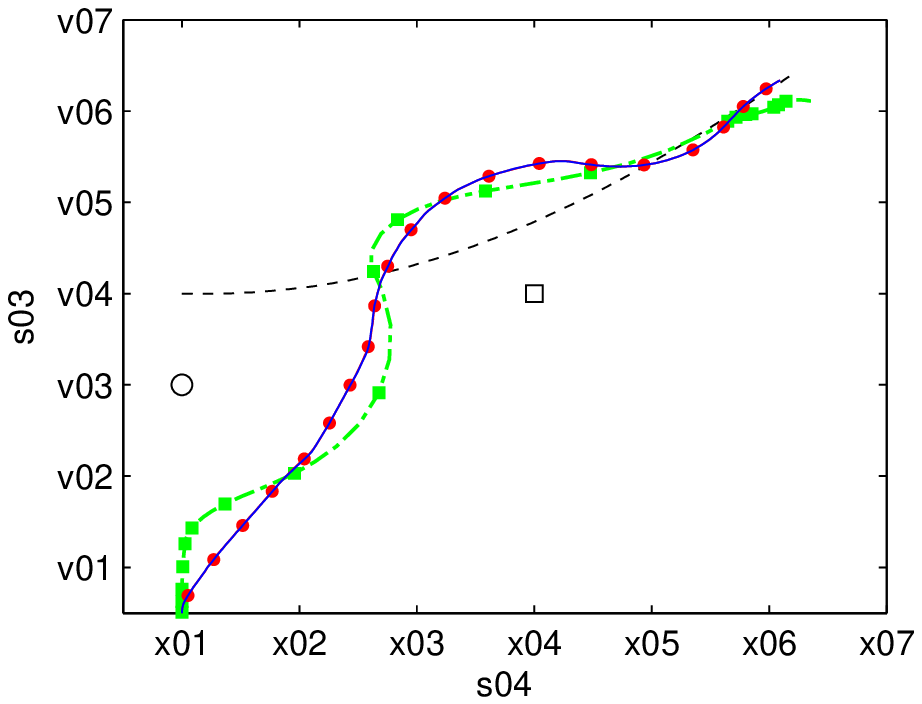}}%
\end{psfrags}%
%
% End f-mobilerobot-trajectory.tex
\caption{Mobile robot trajectory over the simulation:
({\color{Matlabblue}---}) EOC with projection, ({\color{Matlabred}$\bullet$}) EOC no projection,
({\color{Matlabgreen}$\blacksquare$}) Wardi\etal off-line solution~\cite{wardi2010-2137}, ($\cdot$)
target path, ($\circ$) obstacle~1, ($\square$) obstacle 2.}\label{f:mobilerobot-trajectory}\end{figure}

\begin{figure}[]\centering% This file is generated by the MATLAB m-file laprint.m. It can be included
% into LaTeX documents using the packages graphicx, color and psfrag.
% It is accompanied by a postscript file. A sample LaTeX file is:
%    \documentclass{article}\usepackage{graphicx,color,psfrag}
%    \begin{document}\input{f-mobilerobot-eoc-ve-vp}\end{document}
% See http://www.mathworks.de/matlabcentral/fileexchange/loadFile.do?objectId=4638
% for recent versions of laprint.m.
%
% created by:           LaPrint version 3.16 (13.9.2004)
% created on:           18-Feb-2013 13:55:52
% eps bounding box:     10 cm x 7.5 cm
% comment:              
%
\begin{psfrags}%
\psfragscanon%
%
% text strings:
\psfrag{s09}[b][b]{\color[rgb]{0,0,0}\setlength{\tabcolsep}{0pt}\begin{tabular}{c}$\hat{v}_0(t),\tilde{v}_0(t)$\end{tabular}}%
\psfrag{s10}[b][b]{\color[rgb]{0,0,0}\setlength{\tabcolsep}{0pt}\begin{tabular}{c}$\hat{v}_1(t),\tilde{v}_1(t)$\end{tabular}}%
\psfrag{s11}[b][b]{\color[rgb]{0,0,0}\setlength{\tabcolsep}{0pt}\begin{tabular}{c}$\hat{v}_2(t),\tilde{v}_2(t)$\end{tabular}}%
\psfrag{s12}[t][t]{\color[rgb]{0,0,0}\setlength{\tabcolsep}{0pt}\begin{tabular}{c}Time (s)\end{tabular}}%
%
% xticklabels:
\psfrag{x01}[t][t]{0}%
\psfrag{x02}[t][t]{2}%
\psfrag{x03}[t][t]{4}%
\psfrag{x04}[t][t]{6}%
\psfrag{x05}[t][t]{8}%
\psfrag{x06}[t][t]{10}%
\psfrag{x07}[t][t]{12}%
\psfrag{x08}[t][t]{14}%
\psfrag{x09}[t][t]{16}%
\psfrag{x10}[t][t]{18}%
\psfrag{x11}[t][t]{0}%
\psfrag{x12}[t][t]{2}%
\psfrag{x13}[t][t]{4}%
\psfrag{x14}[t][t]{6}%
\psfrag{x15}[t][t]{8}%
\psfrag{x16}[t][t]{10}%
\psfrag{x17}[t][t]{12}%
\psfrag{x18}[t][t]{14}%
\psfrag{x19}[t][t]{16}%
\psfrag{x20}[t][t]{18}%
\psfrag{x21}[t][t]{0}%
\psfrag{x22}[t][t]{2}%
\psfrag{x23}[t][t]{4}%
\psfrag{x24}[t][t]{6}%
\psfrag{x25}[t][t]{8}%
\psfrag{x26}[t][t]{10}%
\psfrag{x27}[t][t]{12}%
\psfrag{x28}[t][t]{14}%
\psfrag{x29}[t][t]{16}%
\psfrag{x30}[t][t]{18}%
%
% yticklabels:
\psfrag{v01}[r][r]{0}%
\psfrag{v02}[r][r]{0.5}%
\psfrag{v03}[r][r]{1}%
\psfrag{v04}[r][r]{0}%
\psfrag{v05}[r][r]{0.5}%
\psfrag{v06}[r][r]{1}%
\psfrag{v07}[r][r]{0}%
\psfrag{v08}[r][r]{0.5}%
\psfrag{v09}[r][r]{1}%
%
% Figure:
\resizebox{8cm}{!}{\includegraphics{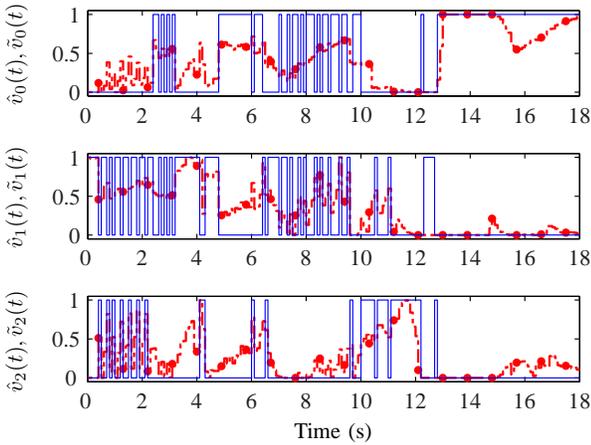}}%
\end{psfrags}%
%
% End f-mobilerobot-eoc-ve-vp.tex
\caption{Mobile robot embedded optimal control mode
selection: ({\color{Matlabblue}---}) projected, ({\color{Matlabred}$\bullet$}) embedded
solution.}\label{f:mobilerobot-eoc-ve-vp}\end{figure}

\Figure~\ref{f:mobilerobot-trajectory} shows the robot trajectories obtained with the embedding approach and using the switching times and mode selection given in~\cite{wardi2010-2137} for the off-line solution.  \Figure~\ref{f:mobilerobot-eoc-ve-vp} shows computed values for the embedding approach and final projected mode switch values. Mode projection is required approximately 87\% of the time.  In Wardi\etal, the mode sequence is predefined as $\{$G2G,Avoid1,G2G,Avoid2,G2G$\}$ and only the switching times are to be found. The embedding approach makes no such assumptions on the mode selection.  The cost obtained with the Wardi\etal offline solution is $22.84$ while the embedding approach here results in a simulation cost of $18.83$, a reduction of 18\%. Further, the embedding approach total solution time was $2.063\dEE{4}$s; Wardi\etal did not report their solution time. As a follow-on test, an MPC approach was implemented with a prediction horizon of $1$ time unit ($10$ partitions of $0.1$ time units). The MPC optimization and implementation resulted in a simulation cost of $20.56$, a total MPC solution time of $71.07$s, and $154$ mode projections (out of 180 partitions); the cost is still about 10\% below the cost of~$22.84$ reported in~\cite{wardi2010-2137}.

%+++++++++++++++++++++++++++++++++++++++++++++++++++++++++++++++++++++++++
\subsection{Two-Tank Hybrid System~\cite{wardi2012acc}}
The two-tank hybrid system in~\cite{wardi2012acc} consists of the draining of one tank into another with the flow rate into the first tank regulated by a valve with two possible values, $\nu=1$ and $\nu=2$. The dynamics are given as
\begin{equation}\label{e:twotank-dyn-ct}
\dot{x}(t)=\bmat{\nu(t)-\sqrt{x_1(t)}\\ \sqrt{x_1(t)}-\sqrt{x_2(t)}}
\end{equation}
where $x=[x_1,x_2]^T$ with $x_1$ the level of tank~1 and $x_2$ the level of tank~2.  The initial conditions are $x_1(0)=x_2(0)=2$ and the states satisfy $x_1,x_2\in[1,4]$. The control problem is to select $\nu(t)$ over $[0,t_f]$, i.e., perform a single optimization, that minimizes
\begin{equation}\label{e:twotank-J}
J(x(0),t_f)=2\int_0^{t_f}\left(x_2(t)-3\right)^2dt
\end{equation}
with $t_f=20$, subject to \Equation\pref{e:twotank-dyn-ct} and the state bounds. The continuous-time dynamics are discretized with the forward-Euler method using a sample interval of $t_s=0.01$. Wardi and Egerstedt~\cite{wardi2012acc} perform the problem optimization with a gradient-descent technique similar to that used for the mobile robot described earlier without pre-defining a mode sequence.  The mode sequence and switch times are calculated to globally minimize the cost function over the entire simulation time.  The approach relies on an insertion gradient that indicates whether or not a change in the mode sequence and/or switch times decreases a cost function. The insertion gradient is iteratively driven toward zero over mode sequences and switch times selected according to their solution algorithm.

Although the PI is suitable for MPP, \MIPT, and CPLEX,  piecewise affine approximations to the square root nonlinearities are needed for approximate solutions using these methods.  Hence the above problem is not solved using MPP, \MIPT, or CPLEX herein. 
% Response to reviewers 1,3: Although the PI is suitable for MPP, \MIPT, and CPLEX,  piecewise affine approximations to the square root nonlinearities are needed for approximate solutions.  Hence the above problem is not solved using MPP, \MIPT, or CPLEX herein.  

The embedding approach only requires the nonlinearities be affine in the continuous control. It is directly applicable to the above optimization using the following embedded model:
\begin{equation}\label{e:twotank-dyn-ct-eocp}
\begin{split}
\dot{\xe}(t)=&(1-\ve(t))\bmat{1-\sqrt{\xe_1(t)}\\ \sqrt{\xe_1(t)}-\sqrt{\xe_2(t)}} \\
&+\ve(t)\bmat{2-\sqrt{\xe_1(t)}\\ \sqrt{\xe_1(t)}-\sqrt{\xe_2(t)}}
\end{split}
\end{equation}
where $\ve\in[0,1]$.  Notice that the embedded problem does not include continuous control inputs, only the mode needs to be selected. The embedded control problem is to minimize \Equation\pref{e:twotank-J} over $\ve$ subject to \Equation\pref{e:twotank-dyn-ct-eocp} and the bounds on the states listed previously.  Appendix~\ref{a:EOCMATSolAlg} outlines how MATLAB's \textit{fmincon} function is used to solve the problem after converting the continuous-time equations to discrete-time equality constraints via the forward-Euler method and approximating the cost function with trapezoidal numerical integration.  Mode projection uses the duty cycle interpretation. For $(1-\ve)t_{s,E}$, mode 0 with $\nu=1$ is applied and then for the remainder of the sample interval mode 1 with $\nu=2$ is used; $t_{s,E}$ is the embedded problem sample interval.
% The embedded problem is solved with a sample interval of $0.25$.

\begin{figure}[]\centering% This file is generated by the MATLAB m-file laprint.m. It can be included
% into LaTeX documents using the packages graphicx, color and psfrag.
% It is accompanied by a postscript file. A sample LaTeX file is:
%    \documentclass{article}\usepackage{graphicx,color,psfrag}
%    \begin{document}\input{f-twotank-states-ve}\end{document}
% See http://www.mathworks.de/matlabcentral/fileexchange/loadFile.do?objectId=4638
% for recent versions of laprint.m.
%
% created by:           LaPrint version 3.16 (13.9.2004)
% created on:           18-Feb-2013 14:07:41
% eps bounding box:     10 cm x 7.5 cm
% comment:              
%
\begin{psfrags}%
\psfragscanon%
%
% text strings:
\psfrag{s09}[b][b]{\color[rgb]{0,0,0}\setlength{\tabcolsep}{0pt}\begin{tabular}{c}$x_1(t)$\end{tabular}}%
\psfrag{s10}[b][b]{\color[rgb]{0,0,0}\setlength{\tabcolsep}{0pt}\begin{tabular}{c}$x_2(t)$\end{tabular}}%
\psfrag{s11}[b][b]{\color[rgb]{0,0,0}\setlength{\tabcolsep}{0pt}\begin{tabular}{c}$\tilde{v}(t)$\end{tabular}}%
\psfrag{s12}[t][t]{\color[rgb]{0,0,0}\setlength{\tabcolsep}{0pt}\begin{tabular}{c}Time (s)\end{tabular}}%
%
% xticklabels:
\psfrag{x01}[t][t]{0}%
\psfrag{x02}[t][t]{2}%
\psfrag{x03}[t][t]{4}%
\psfrag{x04}[t][t]{6}%
\psfrag{x05}[t][t]{8}%
\psfrag{x06}[t][t]{10}%
\psfrag{x07}[t][t]{12}%
\psfrag{x08}[t][t]{14}%
\psfrag{x09}[t][t]{16}%
\psfrag{x10}[t][t]{18}%
\psfrag{x11}[t][t]{20}%
\psfrag{x12}[t][t]{0}%
\psfrag{x13}[t][t]{2}%
\psfrag{x14}[t][t]{4}%
\psfrag{x15}[t][t]{6}%
\psfrag{x16}[t][t]{8}%
\psfrag{x17}[t][t]{10}%
\psfrag{x18}[t][t]{12}%
\psfrag{x19}[t][t]{14}%
\psfrag{x20}[t][t]{16}%
\psfrag{x21}[t][t]{18}%
\psfrag{x22}[t][t]{20}%
\psfrag{x23}[t][t]{0}%
\psfrag{x24}[t][t]{2}%
\psfrag{x25}[t][t]{4}%
\psfrag{x26}[t][t]{6}%
\psfrag{x27}[t][t]{8}%
\psfrag{x28}[t][t]{10}%
\psfrag{x29}[t][t]{12}%
\psfrag{x30}[t][t]{14}%
\psfrag{x31}[t][t]{16}%
\psfrag{x32}[t][t]{18}%
\psfrag{x33}[t][t]{20}%
%
% yticklabels:
\psfrag{v01}[r][r]{0}%
\psfrag{v02}[r][r]{0.5}%
\psfrag{v03}[r][r]{1}%
\psfrag{v04}[r][r]{2}%
\psfrag{v05}[r][r]{3}%
\psfrag{v06}[r][r]{4}%
\psfrag{v07}[r][r]{2}%
\psfrag{v08}[r][r]{3}%
\psfrag{v09}[r][r]{4}%
%
% Figure:
\resizebox{8cm}{!}{\includegraphics{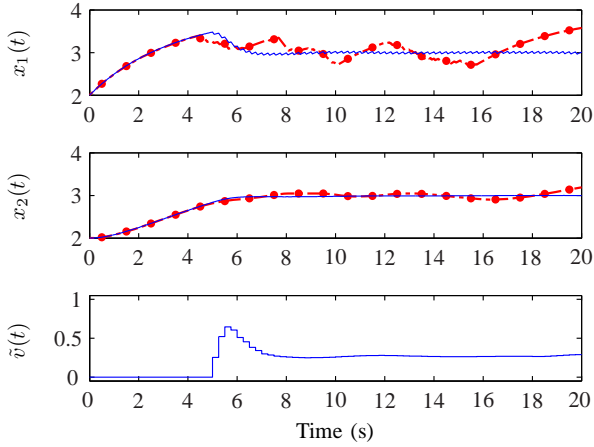}}%
\end{psfrags}%
%
% End f-twotank-states-ve.tex
\caption{Two-tank hybrid system example states and embedded mode switch values: ({\color{Matlabblue}---}) EOC with duty cycle mode projection, ({\color{Matlabred}$\bullet$}) Wardi and Egerstedt solution~\cite{wardi2012acc}.}\label{f:twotank-states-ve}\end{figure}

\Figure~\ref{f:twotank-states-ve} displays the embedding approach performance for $t_{s,E}=0.25$. The $x_2$ state reaches the desired value of $3$ about $6$ time units into the simulation and does not vary more than 1.5\% thereafter. In contrast, the Wardi and Egerstedt results in \Figure~\ref{f:twotank-states-ve} show $x_2$ has variation up to about 7\% from $6$ time units onward and both states show signs of instability as $t$ approaches $t_f$. The reported cost is $4.78$ (recall $t_s=0.01$). On the other hand, for a time step of $0.25$, the embedding approach gives a lower (simulated plant) cost of $4.74$ and shows no signs of instability.  Further, the embedding approach total solution time was $37.56$s; Wardi and Egerstedt did not report their solution time.

%+++++++++++++++++++++++++++++++++++++++++++++++++++++++++++++++++++++++++
\newcommand{\vop}{v'_o}
\newcommand{\ilp}{i'_l}
\subsection{DC-DC Boost Converter Hybrid System~\cite{mariethoz2010-1126}}\label{s:DCDCbooCon}
The DC-DC boost converter~\cite{mariethoz2010-1126} model is a switched lumped parameter circuit, whose continuous-time dynamics on  $[kt_s,(k+1)t_s]$ are
\begin{equation}\label{e:boostconverter-dyn-ct}
\dot{x}(t)=\begin{cases}
    F_1x(t)+f_1v_s,&\; kt_s\leq t \leq (k+d(k))t_s\\
    F_2x(t)+f_2v_s,&\; (k+d(k))t_s<t\leq (k+1)t_s
\end{cases}
\end{equation}
with
\begin{align}
F_1=&H_1\diag\left(-r_l/L,-1/C(r_o+r_c)\right)H_1^{-1}\\
%F_1=&H_1\bmat{-\frac{r_l}{L} & 0 \\ 0 & -\frac{1}{C(r_o+r_c)}}H_1^{-1}\\
f_1=&H_1\bmat{\frac{1}{L} & 0}^T\\
%&H_1=\bmat{ 1 & 0 \\ 0 & \frac{r_o}{r_o+r_c}}\\
&H_1=\diag\left(1,r_o/(r_o+r_c)\right)\\
F_2=&H_2\bmat{-\frac{1}{L}\left(r_l+\frac{r_or_c}{r_o+r_c}\right) & -\frac{r_o}{L(r_o+r_c)} \\
    \frac{r_o}{C(r_o+r_c)} & -\frac{1}{C(r_o+r_c)}}H_2^{-1}\\
f_2=&H_2\bmat{\frac{1}{L} & 0}^T\\
&H_2=\bmat{ 1 & 0 \\ \frac{r_or_c}{r_o+r_c} & \frac{r_o}{r_o+r_c}}
\end{align} where $x=[i_l,v_o]^T$, $i_l$ is the inductor current, $v_o$ the output voltage, $v_s$ the supply voltage, and $d\in[0,0.95]$ the duty cycle in the interval.  The model parameters are $C=100$~$\mu$F, $L=2$~mH, $r_c=0.1$~$\Omega$, $r_l=0.5$~$\Omega$, $r_o=200$~$\Omega$, and the maximum inductor current is ~$2.5$~A; the switching frequency is 20~kHz (the sample interval, $t_s$, is~$50$~$\mu$s).

For the MPC/MPP-based control in~\cite{mariethoz2010-1126}, linearized, discrete-time models are developed from \Equation\pref{e:boostconverter-dyn-ct} corresponding to one of three duty cycle intervals ($i=0,1,2$), $D_0=[0.0.45]$, $D_1=[0.45,0.6]$, and $D_2=[0.6,0.95]$:
\begin{equation}\label{e:boostcovnerter-dyn-linearize-dt}
x'(k+1)=A_{m,i}x'(k)+B_{m,i}d(k)+F_{m,i},
\end{equation}
where $x'(k)=[i'_l,v'_o]^T=[i_l/v_s,v_o/v_s]^T$ (the original state is divided by $v_s$ to avoid the need to generate new linear models if $v_s$ changes~\cite{mariethoz2010-1126}), $x(0)=[0,v_s]^T$, and $v_s=25$~V in this study.  Appendix~\ref{a:bc-dyn-pwa} describes this modeling  process and lists $A_{m,i}$, $B_{m,i}$, and $F_{m,i}$.  The PI~\cite{mariethoz2010-1126} for a horizon window of length~$N$ is
\begin{equation}\label{e:boostconverter-J-eth}
\begin{split}
J(x(k),D_k,N)=&\sum_{k=0}^{N-1}\|Q(v'_o(k+l|k)-v'_{o,ref}(k))\|_p^q \\
    &+\|R(d(k+l|k)-d(k+l-1|k))\|_p^q
\end{split}
\end{equation}
where $D_k=[d(k-1),\ldots,d(k+N-1)]$, $Q$ and $R$ are penalty weight scalars, $v'_{o,ref}(k)=v_{o,ref}/v_s$ is the reference output voltage (assumed constant over $N$), and $p$ designates the p-norm; here $p=1$, $q=1$, and $N\in\{4,6,12\}$.  This MPC problem is also solved with \MIPT, CPLEX, and the embedding approach (using the actual nonlinear model, \Equation\pref{e:boostconverter-dyn-ct}) for a 1-norm PI as in~\cite{mariethoz2010-1126}. However, since the embedding method does not guarantee solution existence or feasibility for p-norm weights on the mode-switch values, $R$ is set to zero whereas $Q=10$. In~\cite{mariethoz2010-1126} the duty cycle is held constant over the $N$ partition prediction horizon. The constant duty cycle over $N$ reduces the optimization to a simple one-dimensional search in contrast to an $N$-dimensional search for a normal MPC horizon as is done herein for all methods.  

Additionally, CPLEX (using \Equation\pref{e:boostcovnerter-dyn-linearize-dt} as specified in Appendix~\ref{a:bc-dyn-pwa}) and the embedding approach (using the actual nonlinear model) are compared using the 1-norm and quadratic PI for the customary $N$ partition MPC window. To use CPLEX for the 1-norm PI, the problem is reformulated as a mixed-integer linear programming problem~\cite{PreConforLinandHybSys2011}. We note that in all the investigations, the application of the duty cycle to the plant is delayed by one sample time interval to be consistent with~\cite{mariethoz2010-1126}.  
% Response to reviewer 3: As the reveiwer correctly points out CPLEX can be used to solve the 1-norm PI problem and we have provided the results in the paper.  Again, thank you for reminding us of the equivalence.  Further, we removed the terminology pseudo-models.

To apply MPP and \MIPT control, $x'(k)$ is augmented with $d(k-1)$ and $v'_{o,ref}$ states to form a new state vector $x''(k)=[i'_l(k),v'_o(k),d(k-1),v'_{o,ref}(k)]^T$; the control input is $\Delta d(k)$. The PI equivalent to \Equation\pref{e:boostconverter-J-eth} is
\begin{equation}\label{e:BooConJ4sta}
%\begin{split}
J(x''(k),\Delta D_k,N)=\sum_{l=0}^{N-1}\|Q'x''(k+l|k)\|_1+\|R\Delta D_k\|_1
%\end{split}
\end{equation}
where $\Delta D_k=[\Delta d(k),\ldots,\Delta d(k+N-1)]^T$, $R=0$, and $Q'(2,2)=10$ and $Q'(2,4)=-10$.  The MPP and \MIPT problem is to minimize \Equation\pref{e:BooConJ4sta} subject to \Equation\pref{e:boostcovnerter-dyn-linearize-dt} and the previously defined bounds on the duty cycle, $i_l$, and $v_o$.

The embedded model is
\begin{equation}\label{e:boostconverter-dyn-ct-eoc}
\begin{split}
\dot{\xe}(t)=&(1-\ve(t))(F_2\xe(t)+f_2v_s) \\
&+\ve(t)(F_1\xe(t)+f_1v_s)
\end{split}
\end{equation}
with PI
\begin{equation}\label{e:boostconverter-J-eoc}
J_E(x(\tpo),\ve,\tpo,\tpf)=\int_{\tpo}^{\tpf}\|Q_E(\ve_o(t)-v_{o,ref}(t))\|_p^qdt
\end{equation}
over $\ve\in[0,1]$ subject to \Equation\pref{e:boostconverter-dyn-ct-eoc} and the previously defined bounds on the duty cycle, $i_l$, and $v_o$.  $Q_E$, $p$, and $q$ are chosen such that the trapezoidal numerical integration approximation of \Equation\pref{e:boostconverter-J-eoc} is consistent with the PI of the MPP, \MIPT, or CPLEX.  As before, the EOCP is solved using MATLAB's \textit{fmincon} function after converting the continuous-time dynamic equations to discrete-time equality constraints via collocation. (See Appendix~\ref{a:EOCMATSolAlg} for the solution method outline.) Mode projection is accomplished using the duty cycle interpretation for each time partition.

\setlength{\tabcolsep}{3.5pt}
\begin{table}[t]
\caption{DC-DC boost converter example EOC (duty cycle interpretation), MPP, \MIPT, and CPLEX 1-norm performance index cost and simulation solution times (``X'' means the algorithm did not converge).}
\label{t:boostconverterSolutionCompare-EMPPMIP}
\centering
\begin{tabular}{l c c c c c}
\hline
\multirow{2}{*}{Method} & \multirow{2}{*}{$N$} & Simulation & Offline & Online & Total\\
 & & Cost & $\tnormSoln$ & $\tnormSoln$ & $\tnormSoln$ \\
\hline
EOC & 4 & 917.33 & - & 1.00 & 1.00 \\ % normalize 85.17
MPP & 4 & 5651.93 & 30.90 & 0.22 & 31.12 \\ % normalize
\MIPT & 4 & 970.09 & 0.08 & 1.17 & 1.26 \\ % normalize
CPLEX & 4 & 946.75& - & 1.15 & 1.15 \\ % normalize 
\hline %& & & & \\
EOC & 6 & 818.45 & - & 9.45 & 9.45 \\ % normalize
MPP & 6 & X & X & X & X\\
\MIPT & 6 & 947.46 & 0.08 & 10.57 & 10.65\\ % normalize
CPLEX & 6 & 946.44 & - & 2.05 & 2.92 \\ % normalize
\hline % & & & & \\
EOC & 12 & 808.94 & - & 44.49 & 44.49 \\ % normalize
MPP & 12 & X & X & X & X\\
\MIPT & 12 & X & X & X & X\\
CPLEX & 12 & 891.79 & - & 1746.65 & 1746.65 \\ % normalize
\hline
\end{tabular}
\end{table}

\begin{table}[t]
\caption{DC-DC boost converter example EOC (duty cycle interpretation) and CPLEX squared 2-norm performance index cost with \Equation\pref{e:boostconverter-J-eoc} and simulation solution times.}
\label{t:boostconverterSolutionCompare-ECpl}
\centering
\begin{tabular}{l c c c}
\hline
Method & $N$ & Cost & $\tnormSoln$ \\
\hline
EOC & 4 & 12593.16  & 1.00\\ % normalize 50.79
CPLEX & 4 & 156547.91 & 1.89 \\ % normalize
\hline % & & & \\
EOC & 6 & 9969.59 & 1.61 \\ % normalize
CPLEX & 6 & 11512.15 & 2.68 \\ % normalize
\hline % & & & \\
EOC & 12 & 9291.08 & 5.26 \\ % normalize
CPLEX & 12 & 10297.26 & 622.28 \\ % normalize
\hline
\end{tabular}
\end{table}

The reference output voltage of
\begin{equation}\label{e:boostconverter-voref}
v_{o,ref}(t)=\begin{cases}
    35~V,& 25~\text{ms}\leq t\leq 35~\text{ms}\\
    50~V,& \text{otherwise}.
\end{cases}
\end{equation}
Longer prediction horizons, such as $N=18$, are not considered since it has been shown that longer prediction horizons do not provide significantly better control~\cite{neely2010thesis}. Table~\ref{t:boostconverterSolutionCompare-EMPPMIP} lists the MPP, \MIPT, CPLEX, and embedding approach 1-norm PI costs, offline solution times (if applicable), control solution times during the simulation, and total control solution times where $\tnormSoln$ is the time value normalized to the least total solution time observed in the table. In several instances, MPP and \MIPT solutions are not listed because control solutions were not generated after 40~hours and the test was stopped.   In all comparable tests, the embedding approach provided the least cost.  Also, as $N$ is increased the cost decreases for all methods having solutions.   The embedding approach also gave the fastest solution times except when compared to CPLEX for $N=6$. We attribute the performance of CPLEX to it being an optimized, commercial  solver.  However, the CPLEX solution time dramatically increases over the embedding approach time when $N=12$, illustrating the ``curse of dimensionality'' associated with mixed-integer programming.  Table~\ref{t:boostconverterSolutionCompare-ECpl} shows the CPLEX and embedding approach quadratic PI costs and total control solution times.  The embedding approach costs are lower than the CPLEX costs in each test.  The embedding approach solution time for $N=4$, $N=6$, and $N=12$ are about 2, 1.6, and 120 times faster than the CPLEX times, respectively.  Boost converter control using the embedding approach implemented in a dedicated microcontroller has been solved as fast as $50$~$\mu$s~\cite{neely2009-1129}.

\begin{figure}[]\centering% This file is generated by the MATLAB m-file laprint.m. It can be included
% into LaTeX documents using the packages graphicx, color and psfrag.
% It is accompanied by a postscript file. A sample LaTeX file is:
%    \documentclass{article}\usepackage{graphicx,color,psfrag}
%    \begin{document}\input{f-boostconverter-eoc-cplex-Hp=12-states-dc}\end{document}
% See http://www.mathworks.de/matlabcentral/fileexchange/loadFile.do?objectId=4638
% for recent versions of laprint.m.
%
% created by:           LaPrint version 3.16 (13.9.2004)
% created on:           18-Feb-2013 14:32:37
% eps bounding box:     10 cm x 7.5 cm
% comment:              
%
\begin{psfrags}%
\psfragscanon%
%
% text strings:
\psfrag{s09}[b][b]{\color[rgb]{0,0,0}\setlength{\tabcolsep}{0pt}\begin{tabular}{c}$v_o(t)$, $v_{o,ref}(t)$\end{tabular}}%
\psfrag{s10}[b][b]{\color[rgb]{0,0,0}\setlength{\tabcolsep}{0pt}\begin{tabular}{c}$i_L(t)$\end{tabular}}%
\psfrag{s11}[b][b]{\color[rgb]{0,0,0}\setlength{\tabcolsep}{0pt}\begin{tabular}{c}Duty Cycle\end{tabular}}%
\psfrag{s12}[t][t]{\color[rgb]{0,0,0}\setlength{\tabcolsep}{0pt}\begin{tabular}{c}Time (ms)\end{tabular}}%
%
% xticklabels:
\psfrag{x01}[t][t]{0}%
\psfrag{x02}[t][t]{5}%
\psfrag{x03}[t][t]{10}%
\psfrag{x04}[t][t]{15}%
\psfrag{x05}[t][t]{20}%
\psfrag{x06}[t][t]{25}%
\psfrag{x07}[t][t]{30}%
\psfrag{x08}[t][t]{35}%
\psfrag{x09}[t][t]{40}%
\psfrag{x10}[t][t]{45}%
\psfrag{x11}[t][t]{50}%
\psfrag{x12}[t][t]{0}%
\psfrag{x13}[t][t]{5}%
\psfrag{x14}[t][t]{10}%
\psfrag{x15}[t][t]{15}%
\psfrag{x16}[t][t]{20}%
\psfrag{x17}[t][t]{25}%
\psfrag{x18}[t][t]{30}%
\psfrag{x19}[t][t]{35}%
\psfrag{x20}[t][t]{40}%
\psfrag{x21}[t][t]{45}%
\psfrag{x22}[t][t]{50}%
\psfrag{x23}[t][t]{0}%
\psfrag{x24}[t][t]{5}%
\psfrag{x25}[t][t]{10}%
\psfrag{x26}[t][t]{15}%
\psfrag{x27}[t][t]{20}%
\psfrag{x28}[t][t]{25}%
\psfrag{x29}[t][t]{30}%
\psfrag{x30}[t][t]{35}%
\psfrag{x31}[t][t]{40}%
\psfrag{x32}[t][t]{45}%
\psfrag{x33}[t][t]{50}%
%
% yticklabels:
\psfrag{v01}[r][r]{0}%
\psfrag{v02}[r][r]{0.5}%
\psfrag{v03}[r][r]{1}%
\psfrag{v04}[r][r]{-5}%
\psfrag{v05}[r][r]{0}%
\psfrag{v06}[r][r]{5}%
\psfrag{v07}[r][r]{20}%
\psfrag{v08}[r][r]{40}%
\psfrag{v09}[r][r]{60}%
%
% Figure:
\resizebox{8cm}{!}{\includegraphics{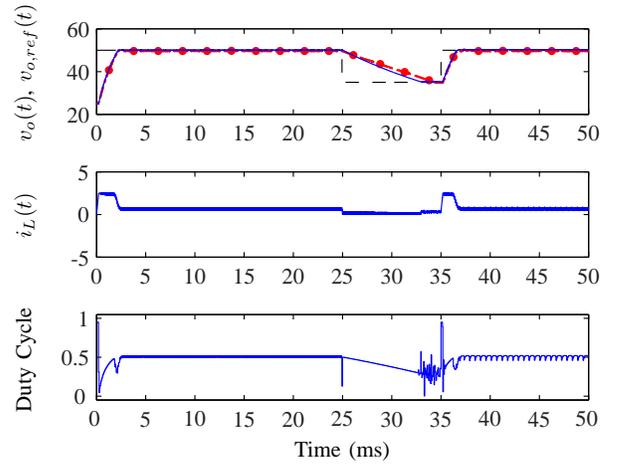}}%
\end{psfrags}%
%
% End f-boostconverter-eoc-cplex-Hp=12-states-dc.tex
\caption{Boost converter hybrid system
states and duty cycle from embedding solution approach with $N=12$: ({\color{Matlabblue}---}) EOC
with duty cycle mode projection, ({\color{Matlabred}$\bullet$}) CPLEX, (--~--)
$v_{o,ref}$.}\label{f:boostconverter-eoc-cplex-Hp=12-states-dc}\end{figure}

\Figure~\ref{f:boostconverter-eoc-cplex-Hp=12-states-dc} shows embedding and CPLEX approach results for $N=12$.  The \Equation\pref{e:boostconverter-voref} reference output voltage appears to be better tracked with the embedding approach during the step change in $v_{o,ref}$ to 35~V on $[25,35)$~ms. The embedding approach is better able to track the 50~V reference on [10,25]~ms (after the initial voltage rise and before the step change in $v_0$) than CPLEX; the average $v_o$ on $[10,25]$~ms from the embedding approach is 50.15~V while that from CPLEX is 49.69~V, the error with the embedding approach is about half that of CPLEX due probably to the approximate model.
%Both control solution approaches are able to track the 50~V output reasonably well between 10 and 25~ms after the initial voltage rise dynamics; the average $v_o$ on $[10,25]$~ms from the embedding approach is 50.15~V while that from CPLEX is 49.69~V.

One final point, good control performance requires an accurate value of the possibly varying load resistance and good estimates of the resistance have been achieved using a nonlinear resistance error system~\cite{neely2010-compel}.  MPP, \MIPT, nor CPLEX has the capability to incorporate a nonlinear estimator directly into the control problem like the embedding approach.

%+++++++++++++++++++++++++++++++++++++++++++++++++++++++++++++++++++++++++
\newcommand{\Xe}{\tilde{X}}
\newcommand{\Ye}{\tilde{Y}}
\newcommand{\thetae}{\tilde{\theta}}
\subsection{Skid-Steered Vehicle Hybrid System~\cite{caldwell2010-5423,caldwell2011-50}}
In~\cite{caldwell2010-5423,caldwell2011-50}, a skid-steered vehicle (SSV) movement hybrid problem is identified as having four modes of operation depending on the lateral sticking or sliding of the wheels: mode~1, front and rear wheels sticking laterally; mode~2, front wheels sticking laterally and back wheels skidding laterally; mode~3, front wheels skidding laterally and back wheels sticking laterally; and mode~4, front and rear wheels skidding laterally.  The equations of motion for modes~1 and~4 are given in~\cite{caldwell2011-50}\footnote{The sign on the $\dot{\theta}$ term in mode~4 is negative here rather than positive as in~\cite{caldwell2011-50} and was changed after correspondence with the authors.}; mode~2 and~3 dynamics are listed in Appendix~\ref{a:ssv-dyn} due to their length.
\begin{align}
f_{1}:=&\begin{cases}
\ddot{X}(t)=\frac{(F_1+F_2+F_3+F_4)\cos\theta(t)}{M}-c_1\dot{X}(t)\\
\ddot{Y}(t)=\frac{(F_1+F_2+F_3+F_4)\sin\theta(t)}{M}-c_1\dot{Y}(t)\\
\ddot{\theta}=0
\end{cases}\\
f_{4}:=&\begin{cases}
  \ddot{X}(t)=\frac{(F_1+F_2+F_3+F_4)\cos\theta(t)}{M}-c_4\dot{X}(t)\\
  \quad+\mu_kg\sin\theta(t)[-\dot{X}(t)\sin\theta(t)+\dot{Y}(t)\cos\theta(t)]\\
  \ddot{Y}(t)=\frac{(F_1+F_2+F_3+F_4)\sin\theta(t)}{M}-c_4\dot{Y}(t)\\
  \quad-\mu_kg\cos\theta(t)[-\dot{X}(t)\sin\theta(t)+\dot{Y}(t)\cos\theta(t)]\\
\ddot{\theta}(t)=\frac{b(F_1-F_2-F_3+F_4)-a^2\mu_KMg\dot{\theta}(t)}{J}
\end{cases}
\end{align}
with $M=m_b+4m_w$,
\begin{equation}\label{e:SkiSteVehIne} %\label{e:SkiSteVehTotalMas}
\begin{split}J=&4m_w(a^2+b^2)+\frac{4m_w}{12}(w_w^2+3w_r^2)\\
    &+\frac{m_b}{12}(B_l^2+B_w^2)
    \end{split}
\end{equation}
where $X$, $Y$, and $\theta$ are position and orientation coordinates in the global coordinate frame; $F_1$/$F_2$ is the right/left rear wheel torque; $F_4$/$F_3$ is the right/left front wheel torque; $c_i$ are mode specific damping coefficients; $\mu_K$ is the coefficient of kinetic friction; $m_b$ is the body mass; $m_w$ is a wheel mass; $B_l$/$B_w$ is the SSV body length/width; $a$/$b$ is half the distance from wheel center to wheel center in the vehicle length/width direction; $w_w$ is half the wheel width; and $w_r$ is the wheel radius.  Model parameters from~\cite{caldwell2010-5423,caldwell2011-50}\footnote{Values for $c_2$, $c_3$, $w_w$, and $w_r$ were obtained through personal correspondence with the authors.} are $c_1=0.7$, $c_2,c_3,c_4=1.2$, $\mu_K=0.8$, $m_b=70$~kg, $m_w=2.5$~kg, $B_l=0.8$~m, $B_w=0.6$~m, $a=0.16$~m, $b=0.28$~m, and $w_w,w_r=0.14$~m.

In~\cite{caldwell2010-5423,caldwell2011-50}, the SSV mode-dependent model is used to generate a reference trajectory, $x_{ref}(t)$, given wheel torques (shown in Figure~3 in~\cite{caldwell2011-50}) and a mode sequence  $(1,4,1)$ with switchings occurring at~5~s and~11~s over a 15~s time interval.  The control problem is to track $x_{ref}(t)$ by minimizing the PI over the controls and modes~\cite{caldwell2011-50}:
\begin{equation}\label{e:ssv-J}
\begin{split}
J(\ve)=&0.5(\xe(t_f)-x_{ref}(t_f))^TP(\xe(t_f)-x_{ref}(t_f))\\
&+\int_{t_0}^{t_f}0.5(\xe(t)-x_{ref}(t))^TQ(\xe(t)-x_{ref}(t))dt
\end{split}
\end{equation}
where $\xe=[\Xe,\dot{\Xe},\Ye,\dot{\Ye},\thetae,\dot{\thetae}]^T$, $Q=\text{diag}(1,1,1,1,1,1)$, $P=\text{diag}(1,10,1,10,1,10)$, $t_0=0$~s, and $t_f=15$~s. 

The methodology for solving this problem coincides with that of the embedding method and is included to illustrate consistency with other results in the literature.  Thus the independently developed methods described in~\cite{caldwell2010-5423,caldwell2011-50} find theoretical justification in~\cite{bengea2005-11}.  Due to the nonlinear nature of the problem, MPP, \MIPT, and CPLEX methods are not utilized without developing piecewise affine approximate models. 

Hence, consistent with the embedding method and~\cite{caldwell2010-5423,caldwell2011-50}, the system dynamics are relaxed:
\begin{equation}\label{e:ssv-dyn-ct}
\dot{\xe}(t)=\sum_{\ml=1}^4\ve_{\ml} f_{\ml}(\xe),\quad \sum_{\ml=1}^4 \ve_{\ml}=1,\quad \ve_{\ml}\in[0,1]\\
\end{equation}

Finally, the active mode is chosen as the mode associated with the greatest value of $\ve_{\ml}$.  This is the mode projection method in Meyer\etal~\cite{meyer2013-1} found to give the least cost in an empirical study of projection methods applied to an example hybrid optimal control problem solved with embedding.

\begin{figure}[]\centering% This file is generated by the MATLAB m-file laprint.m. It can be included
% into LaTeX documents using the packages graphicx, color and psfrag.
% It is accompanied by a postscript file. A sample LaTeX file is:
%    \documentclass{article}\usepackage{graphicx,color,psfrag}
%    \begin{document}\input{f-ssv-eoc-ve-vp}\end{document}
% See http://www.mathworks.de/matlabcentral/fileexchange/loadFile.do?objectId=4638
% for recent versions of laprint.m.
%
% created by:           LaPrint version 3.16 (13.9.2004)
% created on:           30-Nov-2013 20:33:37
% eps bounding box:     10 cm x 7.5 cm
% comment:              
%
\begin{psfrags}%
\psfragscanon%
%
% text strings:
\psfrag{s12}[b][b]{\color[rgb]{0,0,0}\setlength{\tabcolsep}{0pt}\begin{tabular}{c}$v_1(t)$\end{tabular}}%
\psfrag{s13}[b][b]{\color[rgb]{0,0,0}\setlength{\tabcolsep}{0pt}\begin{tabular}{c}$v_2(t)$\end{tabular}}%
\psfrag{s14}[b][b]{\color[rgb]{0,0,0}\setlength{\tabcolsep}{0pt}\begin{tabular}{c}$v_3(t)$\end{tabular}}%
\psfrag{s15}[b][b]{\color[rgb]{0,0,0}\setlength{\tabcolsep}{0pt}\begin{tabular}{c}$v_4(t)$\end{tabular}}%
\psfrag{s16}[t][t]{\color[rgb]{0,0,0}\setlength{\tabcolsep}{0pt}\begin{tabular}{c}Time (s)\end{tabular}}%
%
% xticklabels:
\psfrag{x01}[t][t]{0}%
\psfrag{x02}[t][t]{5}%
\psfrag{x03}[t][t]{10}%
\psfrag{x04}[t][t]{15}%
\psfrag{x05}[t][t]{0}%
\psfrag{x06}[t][t]{5}%
\psfrag{x07}[t][t]{10}%
\psfrag{x08}[t][t]{15}%
\psfrag{x09}[t][t]{0}%
\psfrag{x10}[t][t]{5}%
\psfrag{x11}[t][t]{10}%
\psfrag{x12}[t][t]{15}%
\psfrag{x13}[t][t]{0}%
\psfrag{x14}[t][t]{5}%
\psfrag{x15}[t][t]{10}%
\psfrag{x16}[t][t]{15}%
%
% yticklabels:
\psfrag{v01}[r][r]{0}%
\psfrag{v02}[r][r]{0.5}%
\psfrag{v03}[r][r]{1}%
\psfrag{v04}[r][r]{0}%
\psfrag{v05}[r][r]{0.5}%
\psfrag{v06}[r][r]{1}%
\psfrag{v07}[r][r]{0}%
\psfrag{v08}[r][r]{0.5}%
\psfrag{v09}[r][r]{1}%
\psfrag{v10}[r][r]{0}%
\psfrag{v11}[r][r]{0.5}%
\psfrag{v12}[r][r]{1}%
%
% Figure:
\resizebox{8cm}{!}{\includegraphics{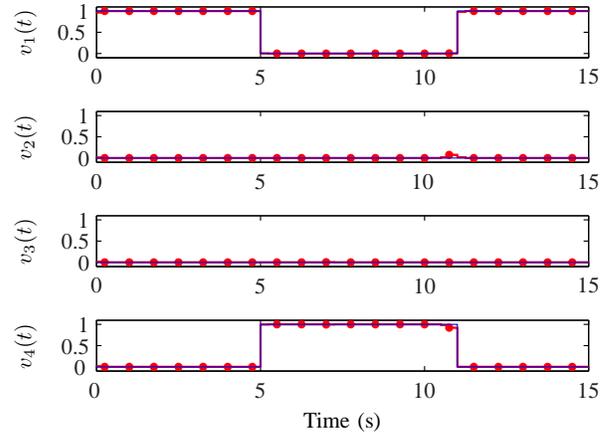}}%
\end{psfrags}%
%
% End f-ssv-eoc-ve-vp.tex
\caption{Skid-steered vehicle embedded solution and projected modes:
({\color{Matlabblue}---}) embedded, ({\color{Matlabred}$\bullet$}) projected.}\label{f:ssv-eoc-ve-vp}\end{figure}

For the embedding approach solution, $[t_0,t_f]$ was divided into 60 equal length partitions. The procedure described in Appendix~\ref{a:EOCMATSolAlg} was then followed with the continuous-time cost computed using trapezoidal numerical integration and the embedded continuous-time dynamics transformed into equality constraints using collocation. The original mode sequence and switching times used to generate the reference state trajectories are recovered as shown in \Figure~\ref{f:ssv-eoc-ve-vp}; mode projection was required for approximately 12\% of the partitions. Also, our embedded solution cost is 2.1\dEE{-4}, meaning the embedded solution tracked the reference trajectories quite well.

%+++++++++++++++++++++++++++++++++++++++++++++++++++++++++++++++++++++++++
\subsection{11 State-Space Region Autonomous Switch Example~\cite{passenberg2010-6666}}\label{s:11StaSpaRegAutSwiExa}
\begin{figure}[]\centering% This file is generated by the MATLAB m-file laprint.m. It can be included
% into LaTeX documents using the packages graphicx, color and psfrag.
% It is accompanied by a postscript file. A sample LaTeX file is:
%    \documentclass{article}\usepackage{graphicx,color,psfrag}
%    \begin{document}\input{f-np-11mode-statespace-Hp=172-up}\end{document}
% See http://www.mathworks.de/matlabcentral/fileexchange/loadFile.do?objectId=4638
% for recent versions of laprint.m.
%
% created by:           LaPrint version 3.16 (13.9.2004)
% created on:           18-Feb-2013 14:55:32
% eps bounding box:     10 cm x 7.5 cm
% comment:              
%
\begin{psfrags}%
\psfragscanon%
%
% text strings:
\psfrag{s03}[l][l]{\color[rgb]{0,0,0}\setlength{\tabcolsep}{0pt}\begin{tabular}{l}1\end{tabular}}%
\psfrag{s04}[l][l]{\color[rgb]{0,0,0}\setlength{\tabcolsep}{0pt}\begin{tabular}{l}2\end{tabular}}%
\psfrag{s05}[l][l]{\color[rgb]{0,0,0}\setlength{\tabcolsep}{0pt}\begin{tabular}{l}3\end{tabular}}%
\psfrag{s06}[l][l]{\color[rgb]{0,0,0}\setlength{\tabcolsep}{0pt}\begin{tabular}{l}4\end{tabular}}%
\psfrag{s07}[l][l]{\color[rgb]{0,0,0}\setlength{\tabcolsep}{0pt}\begin{tabular}{l}5\end{tabular}}%
\psfrag{s08}[l][l]{\color[rgb]{0,0,0}\setlength{\tabcolsep}{0pt}\begin{tabular}{l}6\end{tabular}}%
\psfrag{s09}[l][l]{\color[rgb]{0,0,0}\setlength{\tabcolsep}{0pt}\begin{tabular}{l}7\end{tabular}}%
\psfrag{s10}[l][l]{\color[rgb]{0,0,0}\setlength{\tabcolsep}{0pt}\begin{tabular}{l}8\end{tabular}}%
\psfrag{s11}[l][l]{\color[rgb]{0,0,0}\setlength{\tabcolsep}{0pt}\begin{tabular}{l}9\end{tabular}}%
\psfrag{s12}[l][l]{\color[rgb]{0,0,0}\setlength{\tabcolsep}{0pt}\begin{tabular}{l}10\end{tabular}}%
\psfrag{s13}[l][l]{\color[rgb]{0,0,0}\setlength{\tabcolsep}{0pt}\begin{tabular}{l}11\end{tabular}}%
\psfrag{s14}[b][b]{\color[rgb]{0,0,0}\setlength{\tabcolsep}{0pt}\begin{tabular}{c}$x_2$\end{tabular}}%
\psfrag{s15}[t][t]{\color[rgb]{0,0,0}\setlength{\tabcolsep}{0pt}\begin{tabular}{c}$x_1$\end{tabular}}%
%
% xticklabels:
\psfrag{x01}[t][t]{-10}%
\psfrag{x02}[t][t]{-5}%
\psfrag{x03}[t][t]{0}%
\psfrag{x04}[t][t]{5}%
%
% yticklabels:
\psfrag{v01}[r][r]{-8}%
\psfrag{v02}[r][r]{-6}%
\psfrag{v03}[r][r]{-4}%
\psfrag{v04}[r][r]{-2}%
\psfrag{v05}[r][r]{0}%
\psfrag{v06}[r][r]{2}%
\psfrag{v07}[r][r]{4}%
%
% Figure:
\resizebox{8cm}{!}{\includegraphics{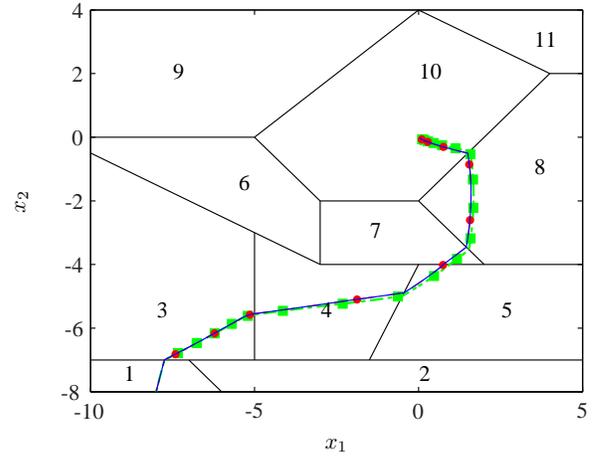}}%
\end{psfrags}%
%
% End f-np-11mode-statespace-Hp=172-up.tex
\caption{Autonomous switched system with 11 a-modes state path
through the labeled state-space partitions: ({\color{Matlabblue}---}) simulation,
({\color{Matlabred}$\bullet$}) states from traditional numerical programming solution with 172
uniform time partitions, ({\color{Matlabgreen}$\blacksquare$}) Passenberg\etal
solution~\cite{passenberg2010-6666}.}\label{f:np-11mode-statespace-Hp=172-up}\end{figure}

Passenberg\etal~\cite{passenberg2010-6666,passenberg2010-4223} proposed a version of a hybrid system minimum principle and applied it to an 11 a-mode linear system.  In the example, a state space, $\Real^2$, is divided into 11 polygonal regions, each associated with a distinct linear state dynamic.  Thus, a discrete state $q(t)\in\{1,2,\ldots,11\}$ identifies a continuous-time state dynamic:
\begin{equation}\label{e:11mode-dyn-ct}
\dot{x}(t)=A_qx(t)+B_qu(t)
\end{equation}
where $x\in\Real^2$, $u\in U_q\subset\Real^2$, $A_q\in\Real^{2\times 2}$ constant, and
$B_q\in\Real^{2 \times 2}$ constant.  The 11 regions are illustrated in
\Figure~\ref{f:np-11mode-statespace-Hp=172-up}.  Each region also has an associated PI integrand:
\begin{equation}
F_q(x,u)=0.5x^T(t)S_q x(t)+0.5u^T(t)R_q u(t)
\end{equation}
where $S_q,R_q\in\Real^{2\times 2}$ are constant penalty weight matrices that depend on $q$.  Given a system trajectory $\left(x(t),u(t)\right)$, the associated discrete state sequence $\{q_1,\ldots,q_\ML\}$ and the switching time sequence $\{t_0,\ldots,t_\ML\}$ ($t_0$ and $t_\ML$ are the initial and final time, respectively), such that $x(t)$ belongs to region $q_\ml$ if $t_{\ml-1}< t < t_\ml$, the PI is then:
\begin{equation}\label{e:11mode-J}
\begin{split}
J(x(t),u(t))=\sum_{\ml=1}^{\ML}\int_{t_{\ml-1}}^{t_\ml } F_{q_\ml} dt.
\end{split}
\end{equation}
The optimal control problem is to drive the states from $x(t_0)=[-8,-8]^T$ at $t_0=0$ to the origin at $t_\ML=2$ such that the PI in \Equation\pref{e:11mode-J} is minimized subject to \Equation\pref{e:11mode-dyn-ct} and constraints on $u$ due to $U_q$ (which is not the same for all values of $q$).  The $A_q$, $B_q$, $S_q$, $R_q$, and $U_q$ are listed in Appendix~\ref{a:11mode-definitions} along with the state-space partition boundary equations. In this 11 state-space regions problem, all the switches are autonomous. At the switching surfaces both the vector fields and the performance index are nonsmooth, which is problematic for traditional optimization methods such as SQP.

The optimization algorithm in~\cite{passenberg2010-4223} can be summarized as follows:
\begin{enumerate}
\item
Choose a feasible sequence of discrete states (adjacent regions) and an appropriate number of switching points (switching time and the value of the state at the switch) on the switching manifolds.\label{Step:Initial}
\item
Find the optimal value for the continuous control $u(t)$ by solving a set of traditional optimization problems for which the switching points provide the fixed initial and final states and times. \label{Step:Optim}
\item
Determine the gradient of the PI at the (autonomous) switching points.\label{Step:Grad}
\item
Based on the gradient computed in Step~\ref{Step:Grad} and the local geometry of the switching manifolds around the switching point, determine the next sequence of the discrete states (regions) and switching points.\label{Step:Update}
\item
Stop if the next switching points are within a certain tolerance of the current ones otherwise return
to Step~\ref{Step:Optim}.
\end{enumerate}
In our view, Step~\ref{Step:Update} is the critical step and the main insight of Passenberg\etal In particular, by separating the computation of the continuous controls from the computation of the switching points, the nonsmoothness at the switching points is completely avoided. Passenberg\etal~\cite{passenberg2010-6666} reported a solution cost of 14.69 for the example with 172 nonuniform time partitions; the solution trajectory is shown in \Figure~\ref{f:np-11mode-statespace-Hp=172-up}\footnote{Passenberg supplied their trajectories in~\cite{passenberg2010-6666} for \Figure~\ref{f:np-11mode-statespace-Hp=172-up} through personal correspondence.}.

In our previous work~\cite{wei2007-264,bengea2011-311} we showed that problems involving a-modes, i.e., autonomous switches, can be solved using traditional mathematical programming methods (in particular, SQP) without any special consideration of the switches and the nonsmoothness associated with them. However, if applied to the example above, the direct application of SQP results in a slightly suboptimal solution.

To evaluate our approach, we first applied traditional numerical programming as suggested in~\cite{wei2007-264} with 172 uniform length partitions;  for better numerical solution behavior, $B_q$ and $R_q$ were scaled so that the continuous controls, $u$, lie in $[-1,1]\times[-1,1]$ for each $q$.  Since all switches are autonomous, no embedding in the vein of~\cite{bengea2011-311} is necessary.  Starting from an initial guess obtained using continuation~\cite{richter1983-660} (which is initialized using Step~\ref{el:fminconinitialguess} in Appendix~\ref{a:EOCMATSolAlg}), we obtained a numerical optimization cost of 14.71 and simulated plant cost with the computed piecewise continuous controls of 14.81, slightly larger than the Passenberg\etal cost of 14.69.

This solution can be improved by accounting for the (discontinuous) transitions across the switching manifolds.  To explicitly account for region transitions, we identified the time partitions in which they occur, estimated the crossover times, appropriately subdivided each identified partition into two of unequal lengths, and then re-optimized. This resulted in a numerical programming cost of 14.69 and simulation cost of 14.71 after 10 iterations on the estimated crossover times, which is consistent with Passenberg\etal reported cost. As an alternative, once the sequence of the discrete states is computed, methods that solve the hybrid optimal control problem assuming a known sequence of discrete states~\cite{wardi2010-2137,xu2000-2683,zefran1997-113} can be used to compute the optimal switching points.

%\begin{table}[t]
%\caption{11-region autonomous switched system problem cost for Passenberg\etal, traditional numerical
%programming (TNP), and TNP with switching point time iteration (TNP/iter).}
%\label{t:11modeSystemCompare-r4concave}
%\centering
%\begin{tabular}{l c c}
%\hline
%Optimization & \multicolumn{2}{c}{Reported Cost} \\
%\hline
%Passenberg\etal\cite{passenberg2010-6666} & \multicolumn{2}{c}{14.69}\\
%& & \\
%\multirow{2}{*}{Optimization} & Optimization & Simulation \\
%& Solution Cost & Cost \\
%\hline
%TNP & 14.71 & 14.81\\
%TNP/iter & 14.69 & 14.71\\
%\hline
%\end{tabular}
%\end{table}

%Table~\ref{t:11modeSystemCompare-r4concave} lists the (similar) costs reported from Passenberg\etal
%and those obtained here.
The trajectory given by Passenberg\etal and the one obtained in the research reported herein using 172 uniform time partitions are plotted in \Figure~\ref{f:np-11mode-statespace-Hp=172-up}; one observes the closeness of the trajectories computed by each of the two algorithms, Passenberg\etal and ours.  The authors would like to acknowledge that the algorithm of Passenberg\etal has an advantage due to the explicit handling of the switching- or discontinuity-manifolds.  In order for classical optimization such as SQP (as asserted in~\cite{wei2007-264}) to obtain similar results, transitions across discontinuity manifolds must be properly accounted for.  Please note MPP, \MIPT, and CPLEX are not applied to this problem because they do not allow for switched PI penalty weights.

%=========================================================================
\section{Discussion and Conclusions}
% This paper delineates specific differences between the embedding
% approach for solving hybrid optimal control problems and
% multi-parametric programming, mixed-integer programming, and
% gradient-descent based methods. First, the embedding approach's
% theoretical underpinnings guarantee the existence of a solution
% under mild conditions. Second, the embedding approach generates a
% convex problem and is solvable with widely available optimization
% packages such as sequential quadratic programming. Third, the
% embedding approach permits nonlinear state models as long as they
% are linear in the continuous time controls. Fourth, the embedding
% approach does not require any assumptions about the mode sequence.

\begin{table*}[t]
\caption{Hybrid optimal control problem solution approaches key differences (GDB is gradient-descent
based methods).}
\label{t:hocpsolutioncompare}
\centering
\begin{tabular}{l c c c c c}
\hline
Characteristic & Embedded & MPP & \MIPT & CPLEX & GDB~\cite{wardi2010-2137}/\cite{wardi2012acc}\\
\hline
General nonlinear control models allowed & Yes & No & No & No & Yes/Yes\\
Switched PI weights allowed & Yes & No & No & No & No/No\\
Utilizes traditional numerical programming & Yes & No & No & No & No/No\\
Computational complexity is NP-hard and rises exponentially with no. of modes & No & Yes & Yes & Yes & No/No\\
Autonomous switches need discrete states assigned & No & Yes & Yes & Yes & -/-\\
Problem initialized with predefined mode sequence & No & No & No & No & Yes/No\\
Real-time control implementation & Yes & Yes & No & No & No/No\\
Shrinking horizon length & Yes & No & Yes & Yes & Yes/Yes\\
Parameter adaptation possible & Yes & No & Yes & Yes & No/No\\
\hline
\end{tabular}
\end{table*}
% Response to reviewer 2: We have changed "Receding horizons" to "Shrinking horizon length".  

The six recently published hybrid optimal control problem examples presented demonstrate the differences and similarities between the embedding approach, MPP, MIP (\MIPT and CPLEX), and gradient based methods as well as the applicability of traditional numerical programming to autonomously switched systems. Table~\ref{t:hocpsolutioncompare} lists key differences between the various SOCP solution approaches.  From the table, the embedding approach addresses a wider class of problems, does not require a specialized solver to implement, does not require additional variables for autonomous switches, and does not require supplying an initial mode sequence.

The embedding approach produced lower performance index costs than all the other optimization methods and found control solutions when other methods failed.  It is noted however that in the presence of autonomous switches this typically requires a finer time partition. The projected embedding produced lower PI costs than the other methods except in two cases: the 25 partition prediction horizon spring-mass and modified-spring-mass examples in which CPLEX produces the lower costs.  The difference here is attributable to two factors: (i) the mode and control projection causing an increase in cost over the embedded cost beyond that of the CPLEX cost, and (ii) the presence of autonomous switches.  In regards to point (ii), an optimal state value (as determined by CPLEX) is at the autonomous switch value and this concurrence is avoided by the numerical solver used by the embedding approach.  This brings up an important consideration.  MPP, \MIPT and CPLEX explicitly parameterize regions of a controlled vector field associated with autonomous switches and hence explicitly deal with the discontinuity that results from an autonomous switch.  In contrast the Jacobian used for the embedding method is discontinuous causing errors in the descent to the optimum.  This suggests that an area of future research is a marriage of the embedding approach with a parameterization of the regions associated with autonomous switches as is done in MIP would result in a superior algorithm with superior convergence properties.  These observations are further born out by the comparison with the 11~autonomous-mode switched system which did not require an embedding. In that comparison, we again found the solution cost obtained with a hybrid minimum principle to be slightly less than that using traditional numerical programming; the important insight is that the numerical algorithm must account for discontinuities in the vector fields either implicitly or explicitly as is done in~\cite{passenberg2010-6666,passenberg2010-4223}. However, we should stress that the basic implementation of the embedding approach as used in this paper produced solutions that were very close to those obtained by these specialized approaches; given its favorable computational complexity and ease of implementation, the embedding approach is thus particularly attractive for engineering applications.

With regards to solution times, the embedding approach gave faster total solution times except in the case of CPLEX for the 3 partition spring-mass and 6 partition boost converter problems.  In regards to real-time implementation, the embedding method must solve the optimization on-line whereas MPP looks up pre-computed feedbacks and critical regions stored in a data base and applies them directly.

Finally, the spring-mass example with $N=25$ demonstrates the need for research to easily construct mode and control projections that better approximates the embedded solution.  The chattering lemma provides a methodology for this to any degree of desired accuracy as is currently being investigated\footnote{Research in this area is currently ongoing at University of California-Berkeley, see arXiv:1208.0062.} to better approximate the embedded solution.

%APPENDICES===============================================================
\appendices

\section{Dynamics Vector Fields}\label{a:DynVecFie}
In regards to the piecewise ${{C}^{1}}$ vector fields, we assume
that the domain ${{\Real}^{n}}\times {{\Real}^{m}}$ is partitioned into
a finite number of nonempty disjoint smooth submanifolds
${{G}_{i}}\,(i=1,\ldots ,{{d}_{G}})$. Each submanifold ${{G}_{i}}$ has
co-dimension ${{p}_{i}}\ge 0$. When ${{p}_{i}}>0$, ${{G}_{i}}$ is
described by a pair of ${{C}^{1}}$ functions
$\sigma_i:\Real^n\to\Real^{p_i}$ and ${{\beta
  }_{i}}:{{\Real}^{n}}\times {{\Real}^{m}}\to {{\Real}^{{{q}_{i}}}}$
(${{q}_{i}}>0$), where the Jacobian $\partial\sigma_{i}/\partial z$
has full rank everywhere. Let
$S_i=\left\{(z,u)\in\Real^n\times\Real^m|\sigma_i(z)=0,\beta_i(z,u_k)<0\right\}$
(the inequalities are taken component-wise), $\partial
S_i=\left\{(z,u)\in\Real^n\times\Real^m|\sigma_i(z)=0,\beta_i(z,u_k)=0\right\}$
and ${{\bar{S}}_{i}}={{S}_{i}}\cup \partial {{S}_{i}}$. We then define
$G_i=S_i\setminus M_{p_i}$ where
$M_k=\underset{p_i>k}{\cup}\bar{S}_i$.  When ${p}_{i}=0$, the
construction of ${{G}_{i}}$ is analogous, except that the function
$\sigma_i$ is absent. In both cases, we define the boundary of
${{G}_{i}}$ to be $\partial G_i=\bar{S}_i\setminus G_i$.  For each
$k\in\{1,\ldots,d_v\}$, the vector field $f_k(x,u_k)$ is $C^1$
on each submanifold $G_i$ and if $(x,{{u}_{k}})\in {{G}_{i}}$ then
$f_k(x,u_k)$ is tangent to ${{G}_{i}}$ at $(x,{{u}_{k}})$. The
trajectory $\left( x(t),{{u}_{k}}(t)\right)$ can only leave $G_i$
through $\partial G_i$. Since discontinuities in
${{f}_{k}}(z,{{u}_{k}})$ determine the autonomous switches and the
boundary points in $\partial {{G}_{i}}$ depend on ${{u}_{k}}$, we need
to impose conditions that prevent instantaneously reversible switches
due to changes in ${{u}_{k}}$. Formally, if for some state value
$x_0\in\Real^n$ there exists values of the control input
${{u}_{0}}$ and ${{u}_{ij}}$ so that $(z_0,u_0)\in G_i$ and
$(z_0,u_{ij})\in\partial G_i$, and at $(z_0,u_{ij})$ a switch occurs
from $G_i$ to $G_j$, then it is not possible to switch from $G_j$ to
$G_i$ at the point $(z_0,\cdot)$ for any value of control $u_k$.

%=========================================================================
\section{MATLAB Embedding Solution Approach Algorithm}\label{a:EOCMATSolAlg}
Here we outline the steps to solve a hybrid optimal control problem with the embedding approach
using MATLAB's \textit{fmincon} function or equivalent.
\begin{enumerate}
\item \label{el:discretizemodel}
(For continuous-time optimal control problem.) Transform the continuous-time embedded system
dynamics into discrete time using the lengths of the prediction horizon partitions and either the
forward-Euler, backward-Euler, or direct
collocation with triangular basis functions~\cite{wei2007-264} methods. The
resulting discretized embedded optimal control problem dynamics are a series of nonlinear equalities.

\item \label{el:discretecost}
(For continuous-time optimal control problem.) Discretize the PI cost over the prediction horizon
partitions using trapezoidal numerical integration.

\item \label{el:discreteconstraints}
Create equality constraints that enforce any supervisory-level
interconnections and embedded mode value sum (embedded mode switch values sum to one) at the
midpoint of each MPC partition if the problem was originally formulated in continuous time or partition
boundary if the problem was given in discrete time.

\item \label{el:fminconinitialguess}
Obtain an initial guess for the \textit{fmincon} variables to be solved for:  states, algebraic
variables, control inputs, and modes.  The initial guess comes from the immediately previous EOCP
solution, i.e., a warm start~\cite{wang2010-267}, when $\tpo>t_0$.  For the first EOCP solution,
$\tpo=t_0$, an in-house created preprocessing function is used to find the initial guess.  This
function works by first solving the EOCP (to a lower numerical tolerance than normal, 1\dEE{-3}
versus 1\dEE{-6}) using \textit{fmincon} for one partition ahead with a user specified initial
guess.  Then, partitions are successively added to the problem, with the EOCP solved after each
addition, until the prediction horizon is reached. During partition addition, the solution
associated with the previously added partition is carried forward to populate the initial guess of
the current EOCP problem.

\item \label{el:fmincon}
Solve the EOCP using \textit{fmincon}.

\item \label{el:vuprojection}
Perform any mode and control projection on the controls for the first partition of the prediction
horizon if the partition's EOCP solution is singular.

\item \label{el:simiulate}
Either obtain a hardware measurement at the end of the current partition
or execute a high accuracy simulation of the continuous-time switched system with the computed
controls over one partition using MATLAB's \textit{ode23t} function.

\item \label{el:feedback}
Slide the starting time of the prediction horizon ahead one partition.  If the end of the
simulation is not reached, return to Step~\ref{el:fminconinitialguess}.
\end{enumerate}

%=========================================================================
\section{Spring-Mass Piece-wise Dynamics}\label{a:SprMasDisTimMod}
The discrete-time matrices for the spring-mass system are in the form
\begin{align}
A_i=\resizematapp{\bmat{a_{11}^i & a_{12}^i \\ a_{21}^i & a_{22}^i}},\;
B_i=\resizematapp{\bmat{b_1^i \\ b_2^i}}\; F_i=\resizematapp{\bmat{f_1^i \\ f_2^i}}. \nonumber
\end{align}
For $i=1$: $a_{11}^1=0.9904$, $a_{12}^1=0.01982$, $a_{21}^1=-a_{12}^1$, $a_{22}^1=-3.965\dEE{-4}$, $b_1^1=9.558\dEE{-3}$, $b_2^1=a_{12}^1$, and $F_1=-B_1$.  For $i=2$: $a_{11}^2=0.9716$, $a_{12}^2=0.01946$, $a_{21}^2=-0.05837$, $a_{22}^2=-1.169\dEE{-3}$, $b_1^2=9.474\dEE{-3}$, $b_2^2=a_{12}^2$, $f_1^2=-0.07105$, and $f_2^2=-0.1459$.  For $i=3$: $a_{11}^3=0.8956$, $a_{12}^3=0.3774$, $a_{21}^3=-a_{12}^3$, $a_{22}^3=0.5183$, $b_1^3=0.1044$, $b_2^3=0.3774$, and $F_3=-B_3$.  For $i=4$: $a_{11}^4=0.6992$, $a_{12}^4=0.3463$, $a_{21}^4=-1.039$, $a_{22}^4=0.3529$, $b_1^4=0.1003$, $b_2^4=a_{12}^4$, $f_1^4=-0.7519$, and $f_2^4=-2.597$.

%\begin{align}
%\begin{split}
%A_{1}=&\resizematapp{\bmat{0.9904 & 0.01982\\-0.01982 & -3.965\dEE{-4}}}\\
%B_{1}=&\resizematapp{\bmat{9.558\dEE{-3}\\0.01982}},\; f_{1}=\resizematapp{\bmat{9.558\dEE{-3}\\0.01982}}
%\end{split}\displaybreak[0]\\
%\begin{split}
%A_{2}=&\resizematapp{\bmat{0.9716 & 0.01946\\-0.05837 & -1.169\dEE{-3}}}\\
%B_{2}=&\resizematapp{\bmat{9.474\dEE{-3}\\0.01946}},\; f_{2}=\resizematapp{\bmat{-0.07105\\-0.1459}}
%\end{split}\displaybreak[0]\\
%\begin{split}
%A_{3}=&\resizematapp{\bmat{0.8956 & 0.3774\\-0.3774 & 0.5183}} \\
%B_{3}=&\resizematapp{\bmat{0.1044\\0.3774}},\;f_{3}=\resizematapp{\bmat{-0.1044\\ -0.3774}}
%\end{split}\displaybreak[0]\\
%\begin{split}
%A_{4}=&\resizematapp{\bmat{0.6992 & 0.3463\\-1.039 & 0.3529}} \\
%B_{4}=&\resizematapp{\bmat{0.1003\\0.3463}},\;f_{4}=\resizematapp{\bmat{-0.7519 \\ -2.5972}}.
%\end{split}
%\end{align}

%=========================================================================
\section{DC-DC Boost Converter Piece-wise Dynamics~\cite{mariethoz2010-1126}}\label{a:bc-dyn-pwa}
The linear DC-DC boost converter dynamics for each duty cycle interval were developed using the
procedure outlined in~\cite{mariethoz2010-1126}.  First, the duty cycle is divided into three
regions: $D_0=[0.0.45]$, $D_1=[0.45,0.6]$, and $D_2=[0.6,0.95]$.  Second, the allowable initial
state conditions are specified. The scaled inductor current, $\ilp=i_l/v_s$, initial condition is
taken to lie in $[0,0.03]$~A/V while the scaled output voltage, $\vop=v_o/v_s$, initial condition
is on $[0,1.6]$. Third, the simulation data for a duty cycle interval is generated. The duty cycle,
inductor current, and output voltage intervals are discretized using five evenly spaced points,
then all possible combinations of the points are simulated over a time interval of $t_s$ duration
using \Equation\pref{e:boostconverter-dyn-ct}. Mariehoz\etal~\cite{mariethoz2010-1126} do not
provide the inductor current or output voltage intervals nor do they give the number of points in
the discretization of the intervals.  Fourth, the simulation data is fit to
\Equation\pref{e:boostcovnerter-dyn-linearize-dt} using least squares. The resulting linear models
associated with the duty cycle intervals have the form
\begin{align}
A_{m,i}=\resizematapp{\bmat{a_{11}^i & a_{12}^i \\ a_{21}^i & a_{22}^i}},\; B_{m,i}=\resizematapp{\bmat{b_1^i \\
b_2^i}}\; F_{m,i}=\resizematapp{\bmat{f_1^i \\ f_2^i}}. \nonumber
\end{align}
For $i=1$: $a_{11}^1=0.9874$, $a_{12}^1=-0.1426$, $a_{21}^1=0.2835$, $a_{22}^1=0.9939$, $b_1^1=0.02828$, $b_2^1=9.680\dEE{-3}$, $f_1^1=0.01347$, and $f_2^1=3.736\dEE{-3}$. For $i=2$: $a_{11}^2=0.9862$, $a_{12}^2=-0.01184$, $a_{21}^2=0.2342$, $a_{22}^2=0.9949$, $b_1^2=0.02000$, $b_2^2=-9.690\dEE{-3}$, $f_1^2=0.01431$, and $f_2^2=0.01074$.  For $i=3$: $a_{11}^3=0.9872$, $a_{12}^3=-5.600\dEE{-3}$, $a_{21}^3=0.1109$, $a_{22}^3=0.9965$, $b_1^3=0.01999$, $b_2^3=-0.01526$, $f_1^3=9.339\dEE{-4}$, and $f_2^3=0.01477$.

%\begin{align}
%\begin{split}
%A_{m,0}=&\resizematapp{\bmat{0.9874 & -0.01426\\0.2835 & 0.9939}},\\
%B_{m,0}=&\resizematapp{\bmat{0.02828\\0.009680}},\; f_{m,0}=\resizematapp{\bmat{0.01347\\3.736\dEE{-3}}}
%\end{split}\displaybreak[0]\\
%\begin{split}
%A_{m,1}=&\resizematapp{\bmat{0.9862 & -0.01184\\0.2342 & 0.9949}}\\
%B_{m,1}=&\resizematapp{\bmat{0.02000\\-0.009690}},\; f_{m,1}=\resizematapp{\bmat{0.01431\\0.01074}}
%\end{split}\displaybreak[0]\\
%\begin{split}
%A_{m,2}=&\resizematapp{\bmat{0.9872 & -5.600\dEE{-3}\\0.1109 & 0.9965}} \\
%B_{m,2}=&\resizematapp{\bmat{0.01999\\-0.01526}},\;\resizematapp{f_{m,2}=\bmat{9.339\dEE{-4} \\ 0.01477}}.
%\end{split}
%\end{align}

%=========================================================================
\section{Skid-Steered Vehicle Hybrid System Dynamics~\cite{caldwell2011-50}}\label{a:ssv-dyn}
\newcommand{\Xdot}{\dot{X}(t)}
\newcommand{\Ydot}{\dot{Y}(t)}
\newcommand{\thetadot}{\dot{\theta}(t)}
\newcommand{\thetat}{\theta(t)}
Skid-steered vehicle mode~2 (front wheels sticking laterally and back wheels skidding laterally)
and~3 (front wheels skidding laterally and back wheels sticking laterally)
dynamics\footnote{Caldwell supplied the mode~2 and~3 dynamics through personal correspondence.}
from~\cite{caldwell2011-50} are given next.  Mode~2 dynamics, $ f_{2}$, are
{\resizeeqnapp\begin{align}
\begin{split}
&\ddot{X}(t)=\frac{1}{M}\left\{\vphantom{\Ydot}(F_1+F_2+F_3+F_4)\cos(\thetat)\right.\\
    &+\frac{\mu_KMg}{2}\sin(\thetat)v_{2,y}\\
    &+\frac{1}{a}\sin(\thetat)[b(F_1-F_2-F_3+F_4)+\frac{\mu_KMga}{2}v_{2,y}]\\
    &-12J\sin(\thetat)[b(F_1-F_2-F_3+F_4)+a^2\mu_KMg\thetadot\\
    &\left.+\alpha_{2}+\beta_{2}]/(a12K_J)\right\} -c_2\Xdot
\end{split}\displaybreak[3]\\
\begin{split}
&\ddot{Y}(t)=[\Xdot+\tan(\thetat)\Ydot]\thetadot\\
    &+\frac{1}{aM}\sin(\thetat)\tan(\thetat)[b(F_1-F_2-F_3+F_4)\\
    &+\frac{a\mu_KMg}{2}v_{2,y}]\\
    &-\frac{\tan(\thetat)}{M}[-(F_1+F_2+F_3+F_4)\cos(\thetat)\\
    &-\frac{\mu_KMg}{2}\sin(\thetat)v_{2,y}] %\\
    -\frac{1}{2aM12K_J}\gamma_{2} %\\
    -c_2\Ydot
\end{split}\displaybreak[3]\\
\begin{split}
&\ddot{\theta}(t)=-12[b(F_1-F_2-F_3+F_4)+a^2\mu_KMg\thetadot \\
&+\alpha_{2}+\beta_{2}]/(12K_J)
\end{split}
\end{align}}
where $M=m_b+4m_w$, $J$ is from \Equation\pref{e:SkiSteVehIne}, {\resizeeqnapp
\begin{align}
\alpha_{2}=&-aM\Xdot(\mu_Kg\sin(\thetat)-\cos(\thetat)\thetadot) \displaybreak[3]\\
\beta_{2}=&aM\Ydot(\mu_Kg\cos(\thetat)+\sin(\thetat)\thetadot) \displaybreak[3]\\
\begin{split}\gamma_{2}=&[12(J+2a^2M)-12 J \cos(2\thetat)]\\
    &\times \sec(\theta(t))[b(F_1-F_2-F_3+F_4)+a^2\mu_KMg\thetadot\\
    &+\alpha_{2}+\beta_{2}]
\end{split}\displaybreak[3]\\
K_J=&J+Ma^2 \displaybreak[3]\\
v_{2,y}=&-\sin(\theta(t))\Xdot+\cos(\theta(t))\Ydot+a\thetadot.
\end{align}}
Mode~3 dynamics, $f_{3}$, are {\resizeeqnapp
\begin{align}
\begin{split}
&\ddot{X}(t)=\frac{1}{M}\left\{\vphantom{\Ydot}(F_1+F_2+F_3+F_4)\cos(\thetat)\right.\\
    &-\frac{\mu_KMg}{2}\sin(\thetat)v_{3,y}\\
    &-\frac{1}{a}\sin(\thetat)[b(F_1-F_2-F_3+F_4)+\frac{a\mu_KMg}{2}v_{3,y}]\\
    &+12J\sin(\thetat)[b(F_1-F_2-F_3+F_4)+a^2\mu_KMg\thetadot\\
    &\left.\alpha_3+\beta_3]/(a12K_J)\right\}-c_3\Xdot
\end{split}\displaybreak[3]\\
\begin{split}
&\ddot{Y}(t)=[\Xdot+\tan(\thetat)\Ydot]\thetadot\\
    &-\frac{1}{(aM)}\sin(\thetat)\tan(\thetat)[b(F_1-F_2-F_3+F_4)\\
    &+\frac{a\mu_KMg}{2}v_{3,y}]\\
    &-\frac{\tan(\thetat)}{M}[-(F_1+F_2+F_3+F_4)\cos(\thetat)\\
    &+\mu_KMg\sin(\thetat)v_{3,y}]
    +\frac{1}{2aM12K_J}\gamma_3
    -c_3\Ydot
\end{split}\displaybreak[3]\\
\begin{split}
&\ddot{\theta}(t)=-12[b(F_1-F_2-F_3+F_4)+a^2\mu_KMg\thetadot\\
    &+\alpha_3+\beta_3]/(12K_J)
\end{split}
\end{align}}
with $\alpha_3=-\alpha_2$, $\beta_3=-\beta_2$, and {\resizeeqnapp
\begin{align}
\begin{split}
\gamma_3=&[12(J+2a^2M)-12J\cos(2\thetat)]\\
    &\times \sec(\thetat)[b(F_1-F_2-F_3+F_4)+a^2\mu_KMg\thetadot \\
    &+\alpha_3+\beta_3]
\end{split} \displaybreak[3]\\
v_{3,y}=&\sin(\theta(t))\Xdot-\cos(\theta(t))\Ydot+a\thetadot.
\end{align}}

%=========================================================================
\section{11-Region Autonomous Switch Example Dynamics and State-Space Divisions~\cite{passenberg2010-6666,passenberg2010-4223}}\label{a:11mode-definitions}

Passenberg\etal~\cite{passenberg2010-6666} 11-region autonomous switched system problem matrices
and state-space partition division lines are listed here\footnote{Passenberg supplied the problem
values through personal correspondence.}. For the $q$-th autonomous mode, $B_q$ is a $2\times 2$
identity matrix (except for $q=10$, then $B_{10}=\smat{1 & 0\\1& 0}$), $S_q=\diag(s_q,s_q)$, $R_q=\diag(r_q,r_q)$ (except for $q=10$,  then $R_{10}=\diag(1,0)$), and
\begin{align}
A_q&=\resizematapp{\bmat{a_{11}^q & a_{12}^q \\ a_{21}^q & a_{22}^q}} \nonumber
\end{align}
The control limits are $U_q\in[-2,2]\times[-2,2]$ for $q\in\{1,\ldots,7,9,11\}$,
$U_8\in[-1,1]\times[-1,1]$, and $U_{10}\in[-10,10]\times[-1,1]$\footnote{$u_2$ dummy variable
defined for consistency with other mode parameter definitions, has no effect on dynamics or cost}.
For $q=1$: $a_{11}^1=-1$, $a_{12}^1=0.5$, $a_{21}^1=-2$, $a_{22}^1=-1$, $s_1=0.6$, $r_1=1$. For
$q=2$: $a_{11}^2=-1$, $a_{12}^2=-2$, $a_{21}^2=0.5$, $a_{22}^2=-0.1$, $s_2=1$, $r_2=1$.  For $q=3$:
$a_{11}^3=-0.5$, $a_{12}^3=-0.1$, $a_{21}^3=-0.1$, $a_{22}^3=-0.1$, $s_3=0.5$, $r_3=0.1$. For
$q=4$: $a_{11}^4=-0.5$, $a_{12}^4=-3$, $a_{21}^4=0.5$, $a_{22}^4=-1$, $s_4=0.5$, $r_4=0.1$. For
$q=5$: $a_{11}^5=-2$, $a_{12}^5=-2$, $a_{21}^5=0.5$, $a_{22}^5=-2$, $s_5=0.1$, $r_5=0.1$. For
$q=6$: $a_{11}^6=-0.1$, $a_{12}^6=0.2$, $a_{21}^6=-2$, $a_{22}^6=-2$, $s_6=8$, $r_6=1$. For $q=6$:
$a_{11}^6=-0.1$, $a_{12}^6=0.2$, $a_{21}^6=-2$, $a_{22}^6=-2$, $s_6=8$, $r_6=1$.  For $q=7$:
$a_{11}^7=-0.5$, $a_{12}^7=-2$, $a_{21}^7=0.5$, $a_{22}^7=-1$, $s_7=1$, $r_7=0.5$.  For $q=8$:
$a_{11}^8=-0.5$, $a_{12}^8=-1$, $a_{21}^8=5$, $a_{22}^8=-0.5$, $s_8=0.1$, $r_8=0.1$.  For $q=9$:
$a_{11}^9=-0.1$, $a_{12}^9=0.5$, $a_{21}^9=-0.5$, $a_{22}^9=-0.2$, $s_9=10$, $r_9=1$.  For $q=10$:
$a_{11}^{10}=-4$, $a_{12}^{10}=2$, $a_{21}^{10}=0$, $a_{22}^{10}=-4$, $s_{10}=2$.  For $q=11$:
$a_{11}^{11}=-1$, $a_{12}^{11}=1$, $a_{21}^{11}=1$, $a_{22}^{11}=-2$, $s_{11}=2$, $r_{11}=1$.

Next, the lines separating the state-space regions are defined. Line~$1$: $x_1=-3$, $x_2\in[-4,-2]$.
Line~$2$: $x_1=-5$, $x_2\in[-7,-3]$. Line~$3$: $x_2=-x_1-14$, $x_1\in[-7,-6]$.  Line~$4$: $x_2=-7$,
$x_1\in[-10,5]$.  Line~$5$: $x_2=-0.5x_1-5.5$, $x_1\in[-10,-3]$.  Line~$6$: $x_2=-4$,
$x_1\in[-3,5]$.  Line~$7$: $x_2=2x_1-4$, $x_1=[-1.5,0]$.  Line~$8$: $x_2=0$, $x_1\in[-10,-5]$.
Line~$9$: $x_2=-x_1-5$, $x_1\in[-5,-3]$.  Line~$10$: $x_2=-2$, $x_1\in[-3,0]$.  Line~$11$:
$x_2=-x_1-2$, $x_1\in[0,2]$.  Line~$12$:  $x_2=x_1-2$, $x_1\in[0,4]$.  Line~$13$: $x_2=2$,
$x_1\in[4,5]$. Line~$14$: $x_2=0.8x_1+4$, $x_1\in[-5,0]$.  Line~$15$: $x_2=-0.5x_1+4$,
$x_1\in[0,4]$.

\section*{Acknowledgments}
This work was partially supported by the Office of Naval Research, and
National Science Foundation grants IIS-0905593, CNS-0910988 and
CNS-1035914. We would like to thank Tim Caldwell and Professor Murphey (skid-steered vehicle), Professor Passenberg (11-region autonomous system), Professor Borrelli (spring-mass), and Professors Egerstedt and Wardi and Philip Twu (mobile robot and two-tank) for supplying invaluable information about the examples.

%We would like to thank Tim Caldwell and Professor Murphey for supplying us with the mode~2 and~3 skid-steered vehicle dynamics
%as well as the vehicle's physical parameters and Professor Passenberg
%for providing their 11-region system problem data and
%results. Professor Borrelli kindly provided access to code and
%information on the spring-mass system. Finally, we would like to thank
%Professors Egerstedt and Wardi as well as Philip Twu who supplied
%essential details of the two-tank and mobile robot problems.

%=========================================================================
\bibliographystyle{IEEEtran}
\bibliography{PhDBibliography}

\end{document}